\documentclass[production,11pt]{article}

\usepackage{amsmath,amsthm,amsfonts,amssymb, latexsym, eepic, floatflt,
  geometry,fancyhdr,cancel,color,picins}
\usepackage[all]{xypic}
\usepackage[dvips]{epsfig}
\newtheorem{thm}{Theorem}[section]
\newtheorem*{claim}{Claim}
\newtheorem{lem}[thm]{Lemma}
\newtheorem{prop}[thm]{Proposition}
\newtheorem{cor}[thm]{Corollary}

\theoremstyle{definition}
\newtheorem*{defin}{Definition}

\DeclareMathOperator{\gr}{gr}
\DeclareMathOperator{\Gr}{Gr}
\DeclareMathOperator{\Sym}{Sym}

\DeclareMathOperator{\ind}{ind}

\def\zz{\mathbb{Z}}
\def\qq{\mathbb{Q}}

\def\torus{\mathbb{T}}

\def\SS{\mathcal{S}}
\def\TT{\mathcal{T}}
\def\bdy{\partial}
\def\tgr{\widetilde{\gr}}

\def\isom{\cong}

\def\spinc{\mathrm{Spin}^c}
\def\spin{\mathrm{Spin}}
\def\SO{\mathrm{SO}}
\def\unitary{\mathrm{U}}
\def\tt{\mathfrak{t}}
\def\ss{\mathfrak{s}}
\def\uu{\mathfrak{u}}
\def\xx{\mathbf{x}}
\def\yy{\mathbf{y}}
\def\ww{\mathbf{w}}
\def\bz{\mathbf{z}}

\newcommand\alphas{\mbox{\boldmath$\alpha$}}
\newcommand\betas{\mbox{\boldmath$\beta$}}
\newcommand\gammas{\mbox{\boldmath$\gamma$}}
\newcommand\deltas{\mbox{\boldmath$\delta$}}

\begin{document}

\title{Covering spaces and $\qq$-gradings on Heegaard Floer homology}
\author{Dan A.\ Lee\\ Duke University\\ dalee@math.duke.edu \and 
Robert Lipshitz\thanks{The second author was supported by a National Science Foundation Graduate Research Fellowship.}
\\ Stanford University\\ lipshitz@math.stanford.edu
}
\date{\today}
\maketitle

\begin{abstract}
Heegaard Floer homology, first introduced by P.~Ozsv\'{a}th and Z.~Szab\'{o} in \cite{holodiskstopo}, associates to a $3$-manifold $Y$ a family of relatively graded abelian groups $HF(Y,\tt)$, indexed by $\spinc$ structures $\tt$ on $Y$.  In the case that $Y$ is a rational homology sphere, Ozsv\'{a}th and Szab\'{o} lift the relative $\zz$-grading to an absolute \mbox{$\qq$-grading} \cite{holotri}. This induces a relative \mbox{$\qq$-grading} on $\bigoplus_{\tt\in\spinc(Y)} HF(Y,\tt)$. In this paper we describe an alternate construction of this relative \mbox{$\qq$-grading} by studying the Heegaard Floer homology of covering spaces.
\end{abstract}

\section{Introduction}
In~\cite{holodiskstopo}, P.~Oszv\'ath and Z.~Szab\'o associated to a $3$-manifold $Y$ families $\widehat{HF}(Y,\tt)$, $HF^+(Y,\tt)$, $HF^-(Y,\tt)$, and $HF^\infty(Y,\tt)$ of abelian groups, indexed by $\spinc$ structures $\tt$ on~$Y$, collectively known as Heegaard Floer homology. (Below, we shall use $HF(Y,\tt)$ to refer to any of these groups.) These groups arise as the homology groups of certain Lagrangian intersection Floer chain complexes $\widehat{CF}(Y,\tt)$, $CF^+(Y,\tt)$, $CF^-(Y,\tt)$, and $CF^\infty(Y,\tt)$, and as such inherit relative $\zz$-gradings or, if $c_1(\tt)$ is non-torsion, relative $\zz/n$ gradings; we denote all of these gradings by $\gr$.

In \cite{holotri}, for $c_1(\tt)$ torsion, Ozsv\'{a}th and Szab\'{o} used a bordism construction to lift $\gr$ to an absolute \mbox{$\qq$-grading} on $HF(Y,\tt)$.  In other words, they found an absolute \mbox{$\qq$-grading} $\tgr$ on $HF(Y,\tt)$ satisfying  $\gr(\xi,\eta)=\widetilde{\gr}(\xi)-\widetilde{\gr}(\eta)$ for all homogeneous elements $\xi,\eta\in HF(Y,\tt)$.  This defines an absolute \mbox{$\qq$-grading} on the group
$$HF(Y,\text{torsion}):=\bigoplus_{c_1(\tt)\text{ is torsion}} HF(Y,\tt).$$
In~\cite{absGrade}, they used this absolute \mbox{$\qq$-grading} to give restrictions on which knots can, under surgery, give rise to lens spaces, and to give restrictions on intersection forms of $4$-manifolds bounding a given $3$-manifold; in~\cite{ManolescuOwens}, C. Manolescu and B. Owens used the absolute \mbox{$\qq$-grading} on the branched double cover of a knot to produce a concordance invariant of knots.

In this paper we use covering spaces to give an alternate construction of the relative \mbox{$\qq$-grading} on $HF(Y,\text{torsion})$ induced by the absolute \mbox{$\qq$-grading} $\widetilde{\gr}$.  Specifically, given homogeneous elements $\xi,\eta\in HF(Y,\text{torsion})$ not necessarily lying in the same $HF(Y,\tt)$, we are able to reconstruct $\widetilde{\gr}(\xi)-\widetilde{\gr}(\eta)$.  See Theorems \ref{main} and \ref{mainagree}.  Of course, the relative \mbox{$\qq$-grading} contains less information about $Y$ than the absolute $\qq$-grading, and therefore our covering space construction is less powerful than the bordism construction of Ozsv\'{a}th and Szab\'{o}; on the other hand, our construction offers a new perspective on $\widetilde{\gr}$ and leads to a simple algorithm for computing the relative \mbox{$\qq$-grading} at the chain level.

Recall (\cite[Section 4]{holodiskstopo}) that there are short exact sequences of chain complexes 
\[0\to\widehat{CF}(Y,\tt)\to CF^+(Y,\tt)\to CF^+(Y,\tt)\to 0
\]
\[
0\to CF^-(Y,\tt)\to CF^\infty(Y,\tt)\to CF^+(Y,\tt)\to 0.
\]
We specify that the relative $\qq$-gradings on these chain complexes be preserved by the maps in these short exact sequences. Consequently, it suffices to define the relative \mbox{$\qq$-grading} on $\widehat{CF}(Y,\text{torsion})$, and we restrict our attention to this chain complex for the rest of the paper.

We now describe the structure of this paper.  In Section 2 we explain our
covering space construction of the relative $\qq$-grading.  This section requires no prior
knowledge of Heegaard Floer homology.  We also explicitly describe how to
compute the relative \mbox{$\qq$-grading} at the chain level.  In Section 3 we briefly review Ozsv\'{a}th and Szab\'{o}'s construction of
$\widetilde{\gr}$, and in Section 4 we prove that our definition of the relative \mbox{$\qq$-grading} agrees with the Ozsv\'{a}th-Szab\'{o} definition.  We defer a necessary computation for lens spaces until Section~5.   Finally, we mention some directions for future research in Section 6.

We thank the referee for a careful reading and many helpful comments.

\section{Gradings and covering spaces}

\subsection{Review of the relative $\zz$-grading}\label{sec:GradingReview}

We begin by defining our use of the term ``grading.''
\begin{defin}
Let $G$ and $S$ be abelian groups.

We say that $f$ is an \emph{absolute $S$-grading on $G$ with homogeneous elements $\mathcal{H}$} if $\mathcal{H}\subset A$  generates $A$ and $f:\mathcal{H}\to S$ such that for each $s\in S$, $f^{-1}(s)\cup\{0\}$ is a subgroup of $A$.

A \emph{relative $S$-grading on $G$ with homogeneous elements $\mathcal{H}$} is an equivalence class of absolute $S$-gradings of $G$ with homogeneous elements $\mathcal{H}$, where two absolute $S$-gradings are equivalent if they differ by a constant in $S$.
\end{defin}
In this paper, $S$ will always be $\zz$ or $\qq$.  If $G$ is free abelian, one can specify an absolute $S$-grading on $G$ with homogeneous elements $\mathcal{H}$ by declaring a basis $\mathcal{U}$ of $G$ to be contained in $\mathcal{H}$ and specifying a function $f:\mathcal{U}\to S$.  One can specify a relative $S$-grading of $G$ by instead specifying a function $F:\mathcal{U}\times\mathcal{U}\to S$ that is \emph{additive} in the sense that $F(x,y)+F(y,z)=F(x,z)$.

As mentioned in the introduction, for any $\spinc$ structure $\tt$
on $Y$ with $c_1(\tt)$ torsion, there is a relative $\zz$-grading, $\gr$, on $\widehat{HF}(Y,\tt)$.  
 We will describe $\gr$ at the
chain level.  That is, given a pointed Heegaard diagram $\mathcal{S}$ for $Y$, we will grade the group $\widehat{CF}(\mathcal{S},\tt)$.  More generally, we will define $\gr$ on $\widehat{CF}(\mathcal{S},\tt)$ when $\mathcal{S}$ is an $\ell$-pointed Heegaard diagram, as defined below.  The material in this subsection has been extracted from Section 3 of~\cite{holoDisksLink}, Section 4 of \cite{L}, and the first half of \cite{holodiskstopo}; we gather it here for the reader's convenience.

\begin{defin}
An \emph{$\ell$-pointed Heegaard diagram} $\mathcal{S}$ is a $4$-tuple $(\Sigma,\alphas,\betas,\mathbf{z})$ where:
\begin{itemize}
\item $\Sigma$ is an oriented surface of genus $g>0$.
\item $\alphas$ is a union of disjoint simple closed curves $\alpha_1,\ldots,\alpha_{g+\ell-1}$ which span a rank $g$ sublattice of $H_1(\Sigma)$.
\item $\betas$ is a union of disjoint simple closed curves $\beta_1,\ldots,\beta_{g+\ell-1}$ which span a rank $g$ sublattice of $H_1(\Sigma)$.
\item $\alphas$ intersects $\betas$ transversely.
\item $\mathbf{z}$ is a collection of points $\{z_1,\ldots,z_\ell\}\subset \Sigma\smallsetminus\alphas\smallsetminus\betas.$
\item Each component of $\Sigma\smallsetminus\alphas$ contains exactly one $z_i$.
\item Each component of $\Sigma\smallsetminus\betas$ contains exactly one $z_i$.
\end{itemize}

The $\alpha$ and $\beta$ circles specify two handlebodies, $U_\alpha$ and $U_\beta$, with boundary $\Sigma$.  We say that $\mathcal{S}$ is an \emph{$\ell$-pointed Heegaard diagram for $Y_{\alpha,\beta}=U_\alpha\cup_\Sigma U_\beta$}.

When $\ell=1$, we say that $\mathcal{S}$ is a \emph{pointed Heegaard diagram}.  (The concept of an $\ell$-pointed Heegaard diagram was introduced in~\cite[Section 3]{holoDisksLink}; earlier papers restricted attention to the $\ell=1$ case. Many theorems that were originally proved in the $\ell=1$ case can be trivially generalized to the case of arbitrary $\ell$.)
\end{defin}

Fix an oriented $3$-manifold $Y$ and a metric on $Y$. Let $f$ be a self-indexing Morse-Smale function on $Y$ with $\ell$ index zero and $\ell$ index three critical points, and choose $\ell$ flowlines of $\nabla f$ connecting the index zero and index three critical points in pairs.  One can then construct an $\ell$-pointed Heegaard diagram $\mathcal{S}$ for $Y$ as follows:\\
\indent$\bullet$ Take $\Sigma=f^{-1}(3/2)$.\\
\indent$\bullet$ Take $\alphas$ to be the intersection of $\Sigma$ with the flowlines leaving index one critical points.\\
\indent$\bullet$ Take $\betas$ to be the intersection of $\Sigma$ with the flowlines entering index two critical points.\\
\indent$\bullet$ Take $\mathbf{z}$ to be the intersection of $\Sigma$ with the $\ell$ flowlines chosen above.\\
In this case we say that $f$ is \emph{compatible} with $\mathcal{S}$.  Observe that, given $\mathcal{S}$, one can always construct a compatible Morse function $f$.

Let $\SS$ be an $\ell$-pointed Heegaard diagram for $Y$.  Define
$\torus_\alpha\cap\torus_\beta$ to be the set of all $(g+\ell-1)$-element subsets
$\xx\subset\alphas\cap\betas$ such that each $\alpha_i$ contains exactly one element
of $\xx$, and each $\beta_i$ contains exactly one element of $\xx$.\footnote{The
  notation comes from thinking of $\torus_\alpha$ as the quotient of the torus
  $\alpha_1\times\cdots\times\alpha_{g+\ell-1}$ by permutations (and similar for
  $\torus_{\beta}$) so that the intersection $\torus_\alpha\cap\torus_\beta$ takes
  place in the $(g+\ell-1)$-fold symmetric product of $\Sigma$.}

We now define a map 
$$s:\torus_\alpha\cap\torus_\beta\to\spinc(Y).$$
Fix a metric on $Y$, and choose a Morse function $f$ compatible with $\SS$.  Then $\xx\in\torus_\alpha\cap\torus_\beta$ determines $g+\ell-1$ flowlines connecting the $g+\ell-1$ index one critical points to the $g+\ell-1$ index two critical points in pairs, and $\mathbf{z}$ determines $\ell$ flowlines connecting the $\ell$ index zero critical points to the $\ell$ index three critical points in pairs.  Consider small neighborhoods of these flowlines containing the critical points.  Notice that $\nabla f$ is nonvanishing outside these neighborhoods, and that, since each neighborhood contains two critical points of opposite parity, we can extend $\nabla f$ to a nonvanishing vector field $V(\xx)$ on all of $Y$.

The vector field $V(\xx)$ reduces the structure group of $TY$ from $\SO(3)$ to
$\SO(2)$, and since $\SO(2)=\unitary(1)\subset \unitary(2)=\spinc(3)$, it follows
that $V(\xx)$ determines a $\spinc$ structure $s(\xx)$ on $Y$.  One can check that
$s(\xx)$ does not depend on the choices of metric, compatible Morse function, and
extension of $\nabla f$.  (Specifically, while different choices will determine a
different vector field $V(\xx)$, the new vector field will be \emph{homologous} to
the original one in the sense of \cite{turaev}. Consequently $s(\xx)$ is
unchanged, as shown in \cite{turaev}.)

Now let $\xx,\yy\in \torus_\alpha\cap\torus_\beta$ where
$\xx=\{x_1,\ldots,x_{g+\ell-1}\}$ and $\yy=\{y_1,\ldots,y_{g+\ell-1}\}$.  Choose $1$-chains $a$ in $\alphas$ and $b$ in $\betas$ such that 
$$\partial
a =\partial b = (y_1+\cdots+y_{g+\ell-1})-(x_1+\cdots+x_{g+\ell-1}).$$
 Then $\partial (a-b)=0$, so $a-b$ descends to an element of $H_1(\Sigma)$, which in turn
defines an element $\epsilon(\xx,\yy)\in H_1(Y)$ via the isomorphism
$$
H_1(Y)\isom
{H_1(\Sigma)\over
[\alpha_1],\ldots,[\alpha_{g+\ell-1}],[\beta_1],\ldots,[\beta_{g+\ell-1}]}.
$$
It is clear that $\epsilon(\xx,\yy)$ is independent of choices of $a$ and $b$.  

There is a nice relationship between $\epsilon$ and $\spinc$ structures. Recall that
the set of (homotopy classes of) $\spinc$ structures on $Y$ forms an affine copy of
$H^2(Y)$.  The following lemma is a generalization of Lemma 2.19 of
\cite{holodiskstopo} to the $\ell$-pointed case, and is closely related to Lemma 3.11
in \cite{holoDisksLink}.
\begin{lem} \label{PD}
For $\xx,\yy\in\torus_\alpha\cap\torus_\beta$, the difference class
$\left(s(\xx)-s(\yy)\right)\in H^2(Y)$ is the Poincar\'{e} dual of  \mbox{$\epsilon(\xx,\yy)\in H_1(Y)$}.
\end{lem}

Let $D_1,\ldots, D_N$ be the closures of the connected components of $\Sigma\smallsetminus\alphas\smallsetminus\betas$, thought of as 2-chains, labeled so that $z_{i}\in D_{i}$ for $i\leq \ell$.  Define $C_2(\mathcal{S})$ to be the group of 2-chains in $\Sigma$ generated by the $D_i$'s, and define $\hat{C}_{2}(\mathcal{S})$ to be the subgroup of $C_2(\mathcal{S})$ generated by the $D_{i}$'s with $i>\ell$, or in other words,  $\hat{C}_{2}(\mathcal{S})$ is generated by the closures of connected components of $\Sigma\smallsetminus\alphas\smallsetminus\betas$ not containing any element of $\mathbf{z}$.

For any $A\in C_2(\mathcal{S})$, define $\partial_\alpha A$ to be the intersection of $\partial A$ with $\alphas$.  For all $\xx,\yy\in
\torus_\alpha\cap\torus_\beta$ where $\xx=\{x_1,\ldots,x_{g+\ell-1}\}$ and
$\yy=\{y_1,\ldots,y_{g+\ell-1}\}$, we define
\begin{eqnarray*} 
\pi_2(\xx,\yy)&=&\left\{A\in C_2(\mathcal{S})\,\left|\, \partial\partial_\alpha A =(y_1+\cdots+y_{g+\ell-1})-(x_1+\cdots+x_{g+\ell-1})\right.\right\}\\
\hat{\pi}_2(\xx,\yy)&=&\pi_2(\xx,\yy)\cap\hat{C}_2(\mathcal{S}).
\end{eqnarray*}

\begin{lem} \label{existsA}
For
$\xx,\yy\in\torus_\alpha\cap\torus_\beta$, if $s(\xx)=s(\yy)$, then
$\hat{\pi}_2(\xx,\yy)$ is nonempty. 
\end{lem} 
\begin{proof} Construct $a$ and
$b$ as in our definition of $\epsilon(\xx,\yy)$.  By the previous lemma, $\epsilon(\xx,\yy)=0$.  By the definition of $\epsilon$, this means that $a-b$ plus some $\alpha$ and $\beta$
circles is zero in $H_1(\Sigma)$ and hence equal to the boundary of some
$C=\sum_{i=1}^{N}c_{i}D_{i}\in C_2(\mathcal{S})$.  Observe that $C\in{\pi}_2(\xx,\yy)$.  For $i\leq\ell$, let $A_{i}$ be the closure of the connected component of $\Sigma\smallsetminus\alphas$ containing $z_{i}$, thought of as a $2$-chain in $C_{2}(\mathcal{S})$.  By the definition of an \mbox{$\ell$-pointed} Heegaard diagram, the coefficient of $D_{i}$ in the expression for $A_{j}$ must be $\delta_{ij}$, so one can check that $C-\sum_{i=1}^{\ell}c_{i}A_{i}\in\hat{\pi}_2(\xx,\yy)$.
\end{proof}

Let  $\mathcal{U}_{\tt}=\{\xx\in
\torus_\alpha\cap\torus_\beta\, |\, s(\xx)=\tt\}$, and define $\widehat{CF}(\mathcal{S},\tt)$ to be the free abelian group generated by $\mathcal{U}_{\tt}$.  We also define $\widehat{CF}(\mathcal{S})=\bigoplus_{\tt\in\spinc(Y)} \widehat{CF}(\mathcal{S},\tt)$.  Our goal in this subsection is to define a relative $\zz$-grading, $\gr$, on
$\widehat{CF}(\mathcal{S},\tt)$ when $\tt$ is torsion.  We declare $\mathcal{U}_{\tt}$ to be
homogeneous so that we just need to define an additive function $\gr:\mathcal{U}_{\tt}\times \mathcal{U}_{\tt}\to\zz$.  Before doing so, we introduce some notation.  For any $A=\sum_{i=1}^Na_iD_i\in C_2(\SS)$, we define the \emph{Euler measure of $A$}, $e(A)$, as follows.   If $D_i$ has $p_i$ (not necessarily distinct) ``vertices'' in  $\alphas\cap\betas$, define $e(D_i)=\chi(D_i)-p_i/4$ where $\chi(D_i)$ is the Euler characteristic. Extend this definition linearly to $C_2(\SS)$, so that
$$e(A)=\sum_{i=1}^N a_i (\chi(D_i)-p_i/4).$$

For any $\xx\in \torus_\alpha\cap\torus_\beta$ and $A\in C_2(\SS)$, we define
$n_{\xx}(A)$ as follows.  Each $x\in\alphas\cap\betas$ ``touches'' four of the
$D_i$'s (with possible repetition); define $n_x(A)$ to be the average of the
coefficients of these four $D_i$'s in the expression for $A$.  Finally, we define
$n_{\xx}(A)=\sum_{x\in\xx}n_x(A)$.

\begin{defin}
Let $\xx,\yy\in \torus_\alpha\cap\torus_\beta$ such that
$s(\xx)=s(\yy)$ is torsion.  By the previous lemma there exists
$A\in \hat{\pi}_2(\xx,\yy)$.  Define $\gr(\xx,\yy)$ by the
formula
\begin{equation}\label{gr}
\gr(\xx,\yy)=e(A)+n_\xx(A)+n_\yy(A).
\end{equation}
\end{defin}
The formula (\ref{gr}), first suggested by Ozsv\'ath and Szab\'o, comes from~\cite[Section 4]{L}, where it was proved to agree with the
standard definition in \cite{holodiskstopo} in terms of the Maslov index.\footnote{The proofs in \cite{L} all deal with the case $\ell=1$ but can easily be extended to the general case using the proof of Theorem \ref{Lpointed}.}   It follows that $\gr$ is additive. For completeness, we include a proof that $\gr$ is well-defined.
\begin{prop}\label{indepA}
In the definition above, $\gr(\xx,\yy)$ does not depend on choice of \mbox{$A\in\hat{\pi}_2(\xx,\yy)$}.
\end{prop}
\begin{proof}
Observe that $A$ is unique up to addition of elements in $\hat{\pi}_2(\xx,\xx)$.  By~\cite[Proposition 7.5]{propertiesandapplications}, for any $P\in\hat{\pi}_2(\xx,\xx)$,
\begin{eqnarray*}
\langle c_{1}(s(\xx)),P\rangle&=&e(P)+2n_{\xx}(P)\\
\langle c_1(s(\yy)),P\rangle&=&e(P)+2n_{\yy}(P).
\end{eqnarray*}
  Since $s(\xx)$ and $s(\yy)$ are torsion, the left sides must vanish and we are left with 
  $$n_{\xx}(P)=n_{\yy}(P)=-{1\over 2}e(P).$$  (Indeed, the reason we deal only with torsion $\spinc$ structures $\tt$ throughout this paper is that we need $\langle c_{1}(\tt),P\rangle$ to be zero.)  Thus, calculating $\gr(\xx,\yy)$ using $A+P$ instead of $A$, we have
\begin{eqnarray*}
\gr(\xx,\yy)&=&e(A+P)+n_\xx(A+P)+n_\yy(A+P)\\
&=&e(A)+e(P)+n_\xx(A)+n_{\xx}(P)+n_\yy(A)+n_{\yy}(P)\\
&=&e(A)+n_\xx(A)+n_\yy(A)
\end{eqnarray*}
\end{proof}

We have now defined a relative $\zz$-grading on $\widehat{CF}(\mathcal{S},\tt)$ when $\tt$ is torsion.  If $\SS$ is weakly admissible in the sense of~\cite[Section 5]{holodiskstopo}, one can make $\widehat{CF}(\mathcal{S})$ into a chain complex whose homology is an invariant of $(Y,\ell)$,  independent of the choice of $\mathcal{S}$ \cite{holoDisksLink}.\footnote{We omit the definition of \emph{weakly admissible} and the definition of the differential because both are peripheral to the content of this paper.  In later sections of this paper, all Heegaard diagrams will be implicitly assumed to be weakly admissible.  Arranging this is never difficult; see \cite[Section 5]{holodiskstopo}.}  When $\ell=1$, we can define the \emph{Heegaard Floer homology group} $$\widehat{HF}(Y):=H(\widehat{CF}(\mathcal{S})).$$
Turning our attention back to the case of general $\ell$, it follows from~\cite[Corollary 4.3]{L} and the definition of the differential that $\gr$ descends to a relative $\zz$-grading on $H(\widehat{CF}(\mathcal{S},\tt))$ when $\tt$ is torsion.

We recall the relationship between the $\ell=1$ case and the general case. The following is~\cite[Theorem 4.5]{holoDisksLink}.
\begin{thm}\label{Lpointed}
Let $\mathcal{S}$ be a weakly admissible  $\ell$-pointed Heegaard diagram for $Y$, and let $\tt\in\spinc(Y)$.  Then
$$H(\widehat{CF}(\mathcal{S},\tt))\isom \widehat{HF}(Y\#(\#^{\ell-1}(S^{1}\times S^{2})),\tt\#\ss_{0}),$$
where $\ss_{0}$ is the unique $\spinc$ structure on $S^{1}\times S^{2}$ with $c_{1}=0$.  

Furthermore, when $\tt$ is torsion, the relative $\zz$-gradings on the two sides are the same.
\end{thm}
\begin{proof}
Given $\mathcal{S}=(\Sigma, \alphas,\betas,\mathbf{z})$, define a pointed Heegaard diagram $\mathcal{S}'=(\Sigma',\alphas',\betas',z')$ as follows.  Define $\Sigma'$ by taking $\Sigma$ and attaching $\ell-1$ tubes connecting $D_{i}\ni z_{i} $ to \mbox{$D_{i+1}\ni z_{i+1}$} for $1\leq i\leq\ell-1$.  Let $\alphas'=\alphas$, $\betas'=\betas$, and $z'=z_{1}$.  Now observe that $\mathcal{S}'$ is a pointed Heegaard diagram for $Y\#(\#^{\ell-1}(S^{1}\times S^{2}))$, and that for each $\xx\in  \torus_\alpha\cap\torus_\beta=\torus_{\alpha'}\cap\torus_{\beta'}$, $s'(\xx)=s(\xx)\#\ss_{0}$, where $s'(\xx)$ is computed with respect to $\mathcal{S}'$ and $s(\xx)$ is computed with respect to $\mathcal{S}$.
Therefore $\widehat{CF}(\mathcal{S})=\widehat{CF}(\mathcal{S}')$, and  $\widehat{CF}(\mathcal{S},\tt)=\widehat{CF}(\mathcal{S}',\tt\#\ss_{0})$ for each \mbox{$\tt\in\spinc(Y)$}.  It is immediate from its definition, which we have omitted, that the differential is the same on both sides.  Finally, it is clear that when $\tt$ is torsion, the relative $\zz$-gradings on $\widehat{CF}(\mathcal{S},\tt)$ and $\widehat{CF}(\mathcal{S}',\tt\#\ss_{0})$ are the same.  The result follows.
\end{proof}

To summarize, we have constructed a map from $\ell$-pointed Heegaard diagrams $\SS$ and torsion $\spinc$ structures $\tt\in\spinc_\mathrm{tor}(Y_{\alpha,\beta})$ to relative $\zz$-gradings on $\widehat{CF}(\mathcal{S},\tt)$,
$$(\mathcal{S},\tt)\mapsto(\widehat{CF}(\mathcal{S},\tt),\gr).$$
Moreover, for each $\ell$, (after restricting to weakly admissible diagrams) this map descends to a map from oriented $3$-manifolds $Y$ with torsion $\spinc$ structures $\tt\in\spinc_\mathrm{tor}(Y)$ to relative $\zz$-gradings on $\widehat{HF}(Y\#(\#^{\ell-1}(S^{1}\times S^{2})),\tt\#\ss_{0})$,
$$(Y,\tt,\ell) \mapsto \left(\widehat{HF}(Y\#(\#^{\ell-1}(S^{1}\times S^{2})),\tt\#\ss_{0}),\gr\right)$$
with the property that $(Y,\tt,\ell)$ and $(Y\#(\#^{\ell-1}(S^{1}\times S^{2})),\tt\#\ss_{0},1)$ produce the same $\gr$.

\subsection{Gradings and covering spaces}
\label{sec:GradingsCoveringSpaces}

Let $p:\tilde{Y}\to Y$ be an $n$-fold, connected covering map.\footnote{All covers of $3$-manifolds will be assumed to be connected, unless explicitly stated otherwise.}   Given an $\ell$-pointed
Heegaard diagram $\mathcal{S}=(\Sigma, \alphas,\betas,\mathbf{z})$ for $Y$ and the covering map
$p$, consider the preimage of $\SS$ under $p$, which we denote by $\tilde{\mathcal{S}}=(\tilde{\Sigma}, \tilde{\alphas},\tilde{\betas},\tilde{\mathbf{z}})$.  We claim that $\tilde{\mathcal{S}}$ is an $n\ell$-pointed Heegaard diagram
for $\tilde{Y}$.

To verify the claim, fix a metric on $Y$ and a Morse function $f$ on $Y$ compatible
with $\SS$.  Pull both of them back to $\tilde{Y}$, so that we have a self-indexing
Morse-Smale function $\tilde{f}$ on $\tilde{Y}$ that has $n\ell$ index zero and
$n\ell$ index three critical points.  Also, $\mathbf{z}$ determines $\ell$ flowlines
for $\nabla f$, which lift to $n\ell$ flowlines for $\tilde{\nabla}\tilde{f}$.  As
described in Section \ref{sec:GradingReview}, $\tilde{f}$ and the $n\ell$ flowlines
determine an  $n\ell$-pointed Heegaard diagram for $\tilde{Y}$. It is clear that this
$n\ell$-pointed Heegaard diagram is exactly $\tilde{\mathcal{S}}$.

We say that $\tilde{\mathcal{S}}$ is a \emph{covering Heegaard diagram} of $\mathcal{S}$.  Note that $\tilde{\Sigma}$ is a connected $n$-fold cover of $\Sigma$ and hence a
surface of genus $ng-n+1$.  Also observe that $\tilde{\alphas}$ is a
disjoint union of $n(g+\ell-1)$ circles, as is $\tilde{\betas}$.

In the situation described above, we have the following useful fact.
\begin{lem}\label{liftspin}
For any $\xx\in \torus_\alpha\cap\torus_\beta$,
$$p^*(s(\xx))=s(\tilde{\xx}),$$
where $\tilde{\xx}$ is the inverse image of $\xx$ under $p$. 
\end{lem}
\begin{proof}
As above, consider a metric $g$ on $Y$, a Morse function $f$ on $Y$ compatible with
$\SS$, and their pullbacks $\tilde{g}$ and $\tilde{f}$ under $p$.  Observe that $\xx$
determines $g+\ell-1$ flowlines of $\nabla f$ connecting the index one and index two
critical points of $f$ in pairs, which must be covered by the $n(g+\ell-1)$ flowlines
of $\tilde{\nabla}\tilde{f}$ determined by $\tilde{\xx}$.  Recalling the definition
of $s$ from Section \ref{sec:GradingReview}, it is now evident that one can choose
the vector fields $V(\xx)$ and $V(\tilde{\xx})$ so that $V(\tilde{\xx})$ is the
pullback of $V(\xx)$. It follows that the diagram
\[
\xymatrix{
 & & B\spinc(3) \ar[d]
\\
\tilde{Y}\ar[r]^{p}\ar@/^1pc/[urr]^(.47){s(\tilde{\xx})}
\ar@/_1pc/[rr]_(.47){T\tilde{Y}} 
& Y\ar[r]^(.4){TY}\ar[ur]^(.47){s(\xx)} & B\SO(3)
}
\]
commutes. This is the desired result.
\end{proof}

The preimage map $p^{-1}:\torus_\alpha\cap\torus_\beta\to   \torus_{\tilde{\alpha}}\cap\torus_{\tilde{\beta}}$ induces a map of groups $\widehat{CF}(\mathcal{S})\to\widehat{CF}(\tilde{\mathcal{S}})$, which we denote $\xi\mapsto\tilde{\xi}$. (Note that in general this is not a chain map.) Let us define
$$\widehat{CF}(\mathcal{S},\text{torsion}):=\bigoplus_{\tt\in\spinc_\mathrm{tor}(Y)}\widehat{CF}(\mathcal{S},\tt).$$
Our goal in this subsection is to prove the following theorem, which completely characterizes our relative $\qq$-grading. 
\begin{thm}\label{main}
There is a unique map from $\ell$-pointed Heegaard diagrams to relative \mbox{$\qq$-gradings} on $\widehat{CF}(\mathcal{S},\mathrm{torsion})$,
$$\mathcal{S}\mapsto(\widehat{CF}(\mathcal{S},\mathrm{torsion}),\Gr)$$
such that the relative \mbox{$\qq$-grading} $\Gr$ has the following properties: 
\begin{itemize}
\item For each $\tt\in\spinc_{\mathrm{tor}}(Y)$, the restriction of $\Gr$ to $\widehat{CF}(\mathcal{S},\tt)$ is equal to $\gr$, where $\gr$ is the relative $\zz$-grading described in the previous subsection.  In particular, they have the same homogeneous elements.
\item If $\tilde{Y}\to Y$ is an $n$-fold cover, then for all homogeneous elements $\xi,\eta \in \widehat{CF}(\mathcal{S},\mathrm{torsion})$,
$$\Gr(\xi,\eta)={1\over n}\Gr(\tilde{\xi},\tilde{\eta}).$$
\end{itemize}

Moreover, for each $\ell$,  (after restricting to weakly admissible diagrams) this map descends to a map from oriented $3$-manifolds $Y$ to relative $\qq$-gradings on $\widehat{HF}(Y\#(\#^{\ell-1}(S^{1}\times S^{2})),\mathrm{torsion})$,
$$(Y,\ell) \mapsto \left(\widehat{HF}(Y\#(\#^{\ell-1}(S^{1}\times S^{2})),\mathrm{torsion}),\Gr\right)$$
with the property that $(Y,\ell)$ and $(Y\#(\#^{\ell-1}(S^{1}\times S^{2})),1)$ produce the same $\Gr$.
\end{thm}

We are now ready to describe our construction of the relative \mbox{$\qq$-grading} on $\widehat{CF}(\mathcal{S},\text{torsion})$.  As in the previous subsection, we simply need to define $\Gr$ on the homogeneous generators in $\torus_\alpha\cap\torus_\beta$.
\begin{defin}
Let $p:\tilde{Y}\to Y$ be an $n$-fold covering map, and let $\xx,\yy\in \torus_\alpha\cap\torus_\beta$ such that $s(\xx)$ and $s(\yy)$ are torsion.  If it happens that $s(\tilde{\xx})=s(\tilde{\yy})$, then we define
$$\Gr(\xx,\yy)={1\over n}\gr(\tilde{\xx},\tilde{\yy}).$$
\end{defin}
In order for this definition to make sense, we need to prove two things.  First, we must show that $\Gr(\xx,\yy)$ is independent of the choice of cover.  Second, we must show that given any $\xx,\yy\in \torus_\alpha\cap\torus_\beta$ with $s(\xx)$ and $s(\yy)$ torsion, there always exists a cover such that $s(\tilde{\xx})=s(\tilde{\yy})$.
\begin{lem}
Let $p:\tilde{Y}\to Y$ be an $n$-fold covering map, and let $\xx,\yy\in \torus_\alpha\cap\torus_\beta$.  If $s(\xx)=s(\yy)$ is torsion, then $s(\tilde{\xx})=s(\tilde{\yy})$ and
$$\gr(\xx,\yy)={1\over n}\gr(\tilde{\xx},\tilde{\yy}).$$
\end{lem}
In other words, computing $\Gr$ using the trivial cover, when possible, is consistent with computing $\Gr$ using any other cover.  
\begin{proof}
The previous lemma shows that $s(\tilde{\xx})=s(\tilde{\yy})$.  To prove the formula, choose \mbox{$A\in \hat{\pi}_2(\xx,\yy)$}.  Now consider the total preimage $\tilde{A}\in {C}_2(\tilde{\mathcal{S}})$ and observe that 
$\tilde{A}\in\hat{\pi}_2(\tilde{\xx},\tilde{\yy})$.  It is clear from the definitions that $e(\tilde{A})=n\cdot e(A)$, $n_{\tilde{\xx}}(\tilde{A})=n\cdot n_{\xx}(A)$, and $n_{\tilde{\yy}}(\tilde{A})=n\cdot n_{\yy}(A)$.  The result now follows from the definition of $\gr$ in equation (\ref{gr}).
\end{proof}

\begin{lem}\label{Lemma:WellDef}
Let $p_1:\tilde{Y}_1\to Y$ be an $n_1$-fold covering map, and let $p_2:\tilde{Y}_2\to Y$ be an $n_2$-fold covering map.  Let $\xx,\yy\in \torus_\alpha\cap\torus_\beta$ such that $s(\xx)$ and $s(\yy)$ are torsion, and suppose that $s(p_1^{-1}\xx)=s(p_1^{-1}\yy)$ and $s(p_2^{-1}\xx)=s(p_2^{-1}\yy)$.  Then 
$${1\over n_1}\gr(p_1^{-1}\xx,p_1^{-1}\yy)={1\over n_2}\gr(p_2^{-1}\xx,p_2^{-1}\yy).$$\end{lem}
\begin{proof}
Let $Y'$ be a connected component of the fibered product $\tilde{Y}_1 \times_Y \tilde{Y}_2$ so that the projections $p_1':Y'\to\tilde{Y}_1$ and $p_2':Y'\to\tilde{Y}_2$ are $n_1'$-fold and $n_2'$-fold covering maps, respectively, and thus $n_1 n_1'=n_2 n_2'$.  Note that since $p_1 p_1'=p_2 p_2'$, given our  $\ell$-pointed Heegaard diagram for~$Y$, we are led to the same covering Heegaard diagram for $Y'$, whether we go through $\tilde{Y}_1$ or~$\tilde{Y}_2$.  Therefore the previous lemma shows that
\begin{eqnarray*}
{1\over n_1}\gr(p_1^{-1}\xx,p_1^{-1}\yy)&=&{1\over n_1 n_1'} 
\gr(p_1'^{-1}p_1^{-1}\xx,p_1'^{-1}p_1^{-1}\yy)\\
&=&{1\over n_2 n_2'}\gr(p_2'^{-1}p_2^{-1}\xx,p_2'^{-1}p_2^{-1}\yy)\\
&=&{1\over n_2}\gr(p_2^{-1}\xx,p_2^{-1}\yy)
\end{eqnarray*}
\end{proof}
Lemma~\ref{Lemma:WellDef} implies that for any $\xx,\yy\in \torus_\alpha\cap\torus_\beta$ with $s(\xx)$ and $s(\yy)$ torsion, $\Gr(\xx,\yy)$ is uniquely defined whenever it is defined at all.  We must now prove that $\Gr(\xx,\yy)$ can always be defined.  For this we will use the following lemma from algebraic topology.
\begin{lem}\label{topology}
Let $Y$ be a connected topological space with the homotopy-type of a CW-complex.  Suppose that $a\in H^2(Y;\zz)$ is $n$-torsion.  Then there exists a $\zz/n$-covering map $p:\tilde{Y}\to Y$ such that $p^*a=0$.
\end{lem}
\begin{proof}
The short exact sequence
\[
\xymatrix{0\ar[r]&\zz\ar[r]^n&\zz\ar[r]&\zz/n\ar[r]&0}
\]
induces the exact sequence
\[
\xymatrix{
H^1(Y;\zz/n)\ar[r]^{{\beta}} & H^2(Y;\zz)\ar[r]^n & H^2(Y;\zz) }.
\]
Since $a$ is $n$-torsion,  $a={\beta}(q)$ for some $q\in H^1(Y;\zz/n)$.  Under the isomorphism $H^1(Y;\zz/n)\isom [Y,K(\zz/n,1)]$, we can represent $q$ by a map $f:Y\to K(\zz/n,1)=B\zz/n$ such that $q=f^*\iota$, where $\iota\in H^1(B\zz/n;\zz/n)$ is the image of the identity via the canonical isomorphism $\text{Hom}(\zz/n,\zz/n)\isom H^1(B\zz/n;\zz/n)$.  
 Let $\pi:E\zz/n\to B\zz/n$ denote the contractible principal $\zz/n$-bundle over $B\zz/n$.  Let $\tilde{Y}$ be the pullback $f^*(E\zz/n)$, so we have the following commutative diagram:
\[
\xymatrix{ \tilde{Y}\ar[r]^(.41)g\ar[d]^p & E\zz/n\ar[d]^\pi\\
Y\ar[r]^(.41)f & B\zz/n
}.\]
We claim that $p:\tilde{Y}\to Y$ is the desired $n$-fold covering map.  Indeed, we have 
$a\in \mathrm{im}({\beta}f^*) =\mathrm{im}(f^*{\beta})$, and thus $p^*a\in\mathrm{im}(p^* f^*) = \mathrm{im} (g^*\pi^*)\subset\mathrm{im} (g^*)=0$, since $H^2(E\zz/n;\zz)=0$.
\end{proof}
\begin{cor}\label{lift}
Let $\xx,\yy\in \torus_\alpha\cap\torus_\beta$ such that $s(\xx)$ and $s(\yy)$ are torsion.  Then for some~$n$, $s(\xx)-s(\yy)$ is $n$-torsion and there is a $\zz/n$-cover $\tilde{Y}\to Y$ such that $s(\tilde{\xx})=s(\tilde{\yy})$.  Consequently, $\Gr(\xx,\yy)$ is well-defined.
\end{cor}
\begin{proof}
Since $c_1(s(\xx))-c_1(s(\yy))=2(s(\xx)-s(\yy))$, the hypotheses imply that $s(\xx)-s(\yy)$ is $n$-torsion for some $n$.  By the previous theorem, there is a $\zz/n$-covering map $p:\tilde{Y}\to Y$ such that $p^*(s(\xx)-s(\yy))=0$.  But $p^*(s(\xx)-s(\yy))=s(\tilde{\xx})-s(\tilde{\yy})$.  
\end{proof}
We have shown that we can always find a cover that allows us to compute $\Gr(\xx,\yy)$, but for a rational homology sphere, there is one cover that always works.  Recall that the maximal abelian cover of $Y$ is the cover corresponding to the commutator subgroup of~$\pi_1(Y)$.
\begin{cor}
Let $Y$ be a rational homology sphere and $\tilde{Y}$ its maximal abelian cover.  Then $\tilde{Y}$ is a finite cover of $Y$, and for all $\xx,\yy\in \torus_\alpha\cap\torus_\beta$, we have $s(\tilde{\xx})=s(\tilde{\yy})$.  Thus $\Gr(\xx,\yy)$ can be computed using this cover.
\end{cor}
\begin{proof}
Since $H_1(Y)\isom\pi_1(Y)/[\pi_1(Y),\pi_1(Y)]$ is finite, $\tilde{Y}$ is a finite cover of $Y$.  The rest of the conclusion follows from the previous corollary, Lemma \ref{liftspin}, and the fact that every cyclic cover is covered by the maximal abelian cover (since every homomorphism from a group $G$ to an abelian group factors through the abelianization of $G$).
\end{proof}

We now turn to additivity.
\begin{cor}
For $\xx,\yy,\mathbf{w}\in \torus_\alpha\cap\torus_\beta$ such that $s(\xx)$, $s(\yy)$, and $s(\mathbf{w})$ are torsion,
$$\Gr(\xx,\yy)+\Gr(\yy,\mathbf{w})=\Gr(\xx,\mathbf{w}).$$
\end{cor}
\begin{proof}
Lemma \ref{topology} and the proof of Corollary \ref{lift} allow us to find a cover where $s(\tilde{\xx})=s(\tilde{\yy})=s(\tilde{\mathbf{w}})$. The additivity of $\gr$ then implies the result.
\end{proof}

This completes the proof of Theorem \ref{main}, except for the last part regarding the invariance of $\Gr$ and the statement that $\Gr$ is essentially independent of $\ell$.  The invariance of $\Gr$ is a direct consequence of the invariance of $\gr$.  The proof that $\Gr$ is independent of $\ell$ follows from the proof that $\gr$ is independent of $\ell$.

\subsection{How to compute}
From the previous subsection, it would appear that every time one wants to compute $\Gr(\xx,\yy)$, one must find an appropriate cover and perform computations in the covering Heegaard diagram $\tilde{\mathcal{S}}$.  However, it turns out that the mere existence of an appropriate cover allows us to perform all of the computations in the original  $\ell$-pointed Heegaard diagram $\mathcal{S}$, without ever thinking about covers at all.

Given $\xx,\yy\in \torus_\alpha\cap\torus_\beta$ such that $s(\xx)$ and $s(\yy)$ are torsion, Corollary \ref{lift} gives us an $n$-fold cover for which we can find $A\in\hat{\pi}_2(\tilde{\xx},\tilde{\yy})$.  Consider the projection $\bar{A}\in {C}_2(\mathcal{S})$, and observe that $\bar{A}\in\hat{\pi}_2(n\xx,n\yy)$ in the sense that  $\bar{A}\in\hat{C}_2(\SS)$ and
$$\partial\partial_{\alpha} \bar{A}=n(y_1+\cdots+y_{g+\ell-1})-n(x_1+\cdots+x_{g+\ell-1}).$$ Also observe that 
\begin{eqnarray*}
\Gr(\xx,\yy)&=&{1\over n}\left[e(A)+n_{\tilde{\xx}}(A)+n_{\tilde{\yy}}(A)\right]
\\
&=&{1\over n}\left[e(\bar{A})+n_\xx(\bar{A})+n_\yy(\bar{A})\right].
\end{eqnarray*}
The moral here is that we do not need to find $A$ in order to compute $\Gr(\xx,\yy)$; it is sufficient to find $\bar{A}$.
\begin{prop}\label{Ntheory}
Let $\xx,\yy\in \torus_\alpha\cap\torus_\beta$ such that $s(\xx)$ and $s(\yy)$ are torsion.  Then $s(\xx)-s(\yy)$ is $n$-torsion for some $n$, and there exists some $A\in\hat{\pi}_2(n\xx,n\yy)$.  For any such $A$,
$$\Gr(\xx,\yy)={1\over n}\left[e(A)+n_\xx(A)+n_\yy(A)\right].$$
\end{prop}
\begin{proof}
The argument above shows that there exists some $A$ for which the result holds.  We just have to prove that the right side of the formula is independent of the choice of $A$.  Observe that $A$ is unique up to addition of elements in $\hat{\pi}_2(\xx,\xx)$.  The rest of the proof is identical to the proof of Proposition \ref{indepA}.
\end{proof}
Therefore computing $\Gr(\xx,\yy)$ is reduced to an elementary exercise in linear algebra:  For each $i >\ell$, compute  $\partial\partial_\alpha D_i$.  Then find integers $a_i$ so that $\sum_{i=\ell+1}^N a_i \partial\partial_\alpha D_i$ equals some multiple $n$ of $(y_1+\cdots+y_{g+\ell-1})-(x_1+\cdots+x_{g+\ell-1})$.  Finally, plug $A=\sum_{i=\ell+1}^{N}a_{i}D_{i}$ into the formula in the previous proposition. 

\section{Bordisms and Ozsv\'{a}th-Szab\'{o}'s construction}
The material in this section is all essentially contained in Ozsv\'ath-Szab\'o's
paper~\cite{holotri}.  However, our purposes require us to articulate more precise chain-level
statements in Section \ref{sec:AbsQGrading}.
\subsection{Review of Heegaard triples}
\begin{defin}An \emph{$\ell$-pointed Heegaard triple} is a 5-tuple $\TT_{\alpha,\beta,\gamma}=(\Sigma, \alphas,\betas,\gammas,\bz)$ such that each of $\SS_{\alpha,\beta}=(\Sigma,\alphas,\betas,\bz)$, $\SS_{\beta,\gamma}=(\Sigma,\betas,\gammas,\bz)$, and $\SS_{\alpha,\gamma}=(\Sigma,\alphas,\gammas,\bz)$ is an $\ell$-pointed Heegaard diagram and $\alphas\cap\betas\cap\gammas=\emptyset$.  Let $Y_{\alpha,\beta}$ be the $3$-manifold specified by $\SS_{\alpha,\beta}$ (and similarly for $Y_{\beta,\gamma}$ $Y_{\alpha,\gamma}$).  When $\ell=1$, $\TT_{\alpha,\beta,\gamma}$ is called a \emph{pointed Heegaard triple.}
\end{defin}

An $\ell$-pointed Heegaard triple specifies a $4$-manifold with boundary $W_{\alpha,\beta,\gamma}$ as follows. Let $T$ denote a triangle. The $\alphas$ circles (respectively $\beta$, $\gamma$) specify a handlebody $U_\alpha$ (respectively $U_\beta$, $U_\gamma$) with boundary $\Sigma$.  Glue each of $U_\alpha\times[0,1]$, $U_\beta\times[0,1]$, and $U_\gamma\times[0,1]$ to an edge of $\Sigma\times T$ (clockwise). The result is a $4$-manifold $W_{\alpha,\beta,\gamma}$ with boundary $-Y_{\alpha,\beta}\cup-Y_{\beta,\gamma}\cup Y_{\alpha,\gamma}$. See~\cite[Section 8.1]{holodiskstopo} for details.

Let $\TT_{\alpha,\beta,\gamma}$ be an $\ell$-pointed Heegaard triple. Let $D_1,\ldots, D_N$ be the closures of the connected components of $\Sigma\smallsetminus\alphas\smallsetminus\betas\smallsetminus\gammas$, thought of as $2$-chains, labeled so that $z_{i}\in D_{i}$ for $i\leq \ell$.  Define $C_2(\TT_{\alpha,\beta,\gamma})$ to be the group of $2$-chains in $\Sigma$ generated by the $D_i$'s, and define $\hat{C}_{2}(\TT_{\alpha,\beta,\gamma})$ to be the subgroup of $C_2(\TT_{\alpha,\beta,\gamma})$ generated by the $D_{i}$'s with $i>\ell$, or in other words,  $\hat{C}_{2}(\TT_{\alpha,\beta,\gamma})$ is generated by the closures of connected components of $\Sigma\smallsetminus\alphas\smallsetminus\betas\smallsetminus\gammas$ not containing any element of $\mathbf{z}$.

For all $\ww\in
\torus_\alpha\cap\torus_\beta$, $\xx\in\torus_\beta\cap\torus_\gamma$, and $\yy\in\torus_\alpha\cap\torus_\gamma$, define
\begin{eqnarray*} 
\pi_2(\ww,\xx,\yy)&=&
\left\{B\in C_2(\TT_{\alpha,\beta,\gamma})\,\left|\, 
\begin{array}{l}
\partial\partial_\alpha B =(y_1+\cdots +y_{g+\ell-1})- (w_1+\cdots+ w_{g+\ell-1})\\  
\partial\partial_\beta B=(w_1+\cdots+ w_{g+\ell-1})-(x_1+\cdots+ x_{g+\ell-1})
\end{array}\right.\right\}\\
\hat{\pi}_2(\ww,\xx,\yy)&=&\pi_2(\ww,\xx,\yy)\cap\hat{C}_2(\TT_{\alpha,\beta,\gamma}).
\end{eqnarray*}
Observe that addition gives a map 
$$\pi_2(\ww',\ww)\times\pi_2(\xx',\xx)\times\pi_2(\yy',\yy)\times\pi_2(\ww,\xx,\yy)\to \pi_2(\ww',\xx',\yy').$$

For $B\in{\pi}_2(\ww,\xx,\yy)$ and $B'\in{\pi}_2(\ww',\xx',\yy')$, say that $B$ and $B'$ are \emph{$\spinc$-equivalent} if there are $A_{\alpha,\beta}\in{\pi}_2(\ww',\ww)$, $A_{\beta,\gamma}\in{\pi}_2(\xx',\xx)$ and $A_{\alpha,\gamma}\in{\pi}_2(\yy,\yy')$ such that $B'=A_{\alpha,\beta}+A_{\beta,\gamma}+A_{\alpha,\gamma}+B$.

\begin{lem}
There is a map $s:\pi_2(\ww,\xx,\yy)\to \spinc(W_{\alpha,\beta,\gamma})$ with the following properties:
\begin{itemize}
\item $\ss\in\mathrm{im}(s)$ if and only if $\ss|_{Y_{\alpha,\beta}}=s(\ww)$, $\ss|_{Y_{\beta,\gamma}}=s(\xx)$, and $\ss|_{Y_{\alpha,\gamma}}=s(\yy)$.
\item For $B\in\pi_2(\ww,\xx,\yy)$ and $B'\in\pi_2(\ww',\xx',\yy')$,
$s(B)=s(B')$ if and only if $B$ and $B'$ are $\spinc$-equivalent.
\end{itemize}
\end{lem}
For a definition of the map $s$ and the proof of this lemma, see \cite[Section 8.1]{holodiskstopo}. 

Let $\hat{\pi}_2^\ss(\ww,\xx,\yy)=\{B\in\hat{\pi}_2(\ww,\xx,\yy)\,|\,s(B)=\ss\}$.
For any $\ss\in\spinc(W_{\alpha,\beta,\gamma})$, we say that an $\ell$-pointed Heegaard triple $\TT_{\alpha,\beta,\gamma}$ \emph{realizes} $\ss$ when there exist $\ww$, $\xx$, and $\yy$ such that $\hat{\pi}_2^\ss(\ww,\xx,\yy)\neq\emptyset$.  The
previous lemma
 implies that $\hat{\pi}_2^\ss(\ww,\xx,\yy)\neq\emptyset$ if and only if $\ss|_{Y_{\alpha,\beta}}=s(\ww)$, $\ss|_{Y_{\beta,\gamma}}=s(\xx)$, and $\ss|_{Y_{\alpha,\gamma}}=s(\yy)$.

 Given an element $B\in\pi_2(\ww,\xx,\yy)$, one can assign an index $\ind(B)\in\zz$
 to $B$; this is the index of the $\overline{\partial}$-operator on some space of
 holomorphic curves.\footnote{For the definition of $\ind$, see \cite[Section 8]{holodiskstopo}, where the notation $\mu$
 is used instead of $\ind$.} This index has the property that for
 $B\in\pi_2(\ww,\xx,\yy)$, $A_{\alpha,\beta}\in\pi_2(\ww',\ww)$,
 $A_{\beta,\gamma}\in\pi_2(\xx',\xx)$ and $A_{\alpha,\gamma}\in\pi_2(\yy',\yy)$,
\[
\ind(A_{\alpha,\beta}+A_{\beta,\gamma}+A_{\alpha,\gamma}+B)=\ind(A_{\alpha,\beta})+\ind(A_{\beta,\gamma})+\ind(A_{\alpha,\gamma}) +\ind(B).
\]
where $\ind(A_{\alpha,\beta})=e(A_{\alpha,\beta})+n_{\ww'}(A_{\alpha,\beta})+n_{\ww}(A_{\alpha,\beta})$, and similarly for $\ind(A_{\beta,\gamma})$ and $\ind(A_{\alpha,\gamma})$.

\begin{defin}
Let $\ww\in
\torus_\alpha\cap\torus_\beta$, $\xx\in\torus_\beta\cap\torus_\gamma$, and $\yy\in\torus_\alpha\cap\torus_\gamma$ such that $s(\ww)$, $s(\xx)$, and $s(\yy)$ are torsion, and let $\ss\in\spinc(W_{\alpha,\beta,\gamma})$ such that $\hat{\pi}_2^\ss(\ww,\xx,\yy)\neq\emptyset$.  Then for any $B\in\hat{\pi}^\ss_2(\ww,\xx,\yy)$, define
$$\gr^\ss(\ww,\xx,\yy)=\ind(B).$$
\end{defin}
The definition is independent of choice of $B$ by the same reasoning as in Proposition~\ref{indepA}, but note that $\gr^\ss(\ww,\xx,\yy)$ is only defined when $\hat{\pi}_2^\ss(\ww,\xx,\yy)\neq\emptyset$.  This ``relative grading'' $\gr$ is additive in the sense that 
$$\gr^\ss(\ww',\xx',\yy')=\gr(\ww',\ww) +\gr(\xx',\xx) + \gr(\yy,\yy') +    \gr^\ss(\ww,\xx,\yy).$$

In order to prove that $\Gr$ agrees with the relative grading induced by
$\widetilde{\gr}$, we will need to use covering Heegaard triples; we collect a few
basic facts about these here.  Given an $\ell$-pointed Heegaard triple
$\TT_{\alpha,\beta,\gamma}=(\Sigma,\alphas,\betas,\gammas,\bz)$ and an $n$-fold
covering $\pi:\tilde{W}\to W_{\alpha,\beta,\gamma}$, we can define a new
$n\ell$-pointed Heegaard triple
$\tilde{\TT}_{\tilde{\alpha},\tilde{\beta},\tilde{\gamma}}=(\tilde{\Sigma},\tilde{\alphas},\tilde{\betas},\tilde{\gammas},\tilde{\bz})$
for $\tilde{W}$ by taking $\tilde{\TT}_{\tilde{\alpha},\tilde{\beta},\tilde{\gamma}}$
to be the preimage of $\TT_{\alpha,\beta,\gamma}$ under $\pi$, as we did in
Section~\ref{sec:GradingsCoveringSpaces} for ordinary Heegaard diagrams. That
$\tilde{\TT}_{\tilde{\alpha},\tilde{\beta},\tilde{\gamma}}$ is a Heegaard triple for
$\tilde{W}$ follows by exactly the same argument as in
Section~\ref{sec:GradingsCoveringSpaces}. We say that
$\tilde{\TT}_{\tilde{\alpha},\tilde{\beta},\tilde{\gamma}}$ is a \emph{covering
  Heegaard triple} of $\TT_{\alpha,\beta,\gamma}$.

\begin{lem} \label{Lemma:CovTrip}Fix $\ww\in \torus_{\alpha}\cap\torus_{\beta}$, $\xx\in\torus_{\beta}\cap\torus_{\gamma}$, $\yy\in\torus_{\beta}\cap\torus_{\gamma}$, and $\ss\in\spinc(W_{\alpha,\beta,\gamma})$.  For any $B\in\hat{\pi}^\ss_2(\ww,\xx,\yy)$ and its total preimage $\tilde{B}$, 
\begin{enumerate}
\item \label{item:s}$s(\tilde{B})=\pi^*s(B)$.
\item \label{item:ind}$\ind(\tilde{B})=n\cdot \ind(B)$.
\end{enumerate}
In particular, if the restriction of $\ss$ to $\partial W_{\alpha,\beta,\gamma}$ is torsion and $\hat{\pi}_2^\ss(\ww,\xx,\yy)\neq\emptyset$, then $\gr^{\pi^*\ss}(\tilde{\ww},\tilde{\xx},\tilde{\yy})=n\cdot \gr^\ss(\ww,\xx,\yy)$.
\end{lem}
\begin{proof}
  The proof of statement~\ref{item:s} follows from the definition of $s(B)$ (which
  we have not given) in an exactly analogous way to Lemma~\ref{liftspin}.

  To prove statement~\ref{item:ind}, recall from~\cite{holodiskstopo} that $\ind$
  denotes the index of the
  linearized $\overline{\partial}$-operator at an appropriate map $\phi$ of a
  triangle $\Delta$ to $\Sym^{g+\ell-1}(\Sigma)$. There is a map
  $u:S\to\Sigma\times\Delta$ for some Riemann surface $S$ which \emph{tautologically
    corresponds} to $\phi$ (see, \textit{e.g.}, \cite[Section 13]{L}).  The index of the
  $\overline{\partial}$-operators at $u$ and $\phi$ agree. (Indeed, there is even an
  identification of index bundles -- see~\cite[Section 13]{L} -- though we do not
  need this stronger result.)

  In our situation, the covering map $\tilde{\Sigma}\to\Sigma$ induces an
  \emph{inclusion} $\Sym^{g+\ell-1}(\Sigma)\to\Sym^{ng+n\ell-n}(\tilde{\Sigma})$ and
  a covering map $\tilde{\Sigma}\times\Delta\to\Sigma\times\Delta$. The maps $\phi$
  and $u$ then induce maps $\tilde{\phi}:\Delta\to\Sym^{ng+n\ell-n}(\tilde{\Sigma})$ and
  $\tilde{u}:\tilde{S}\to\tilde{\Sigma}\times\Delta$. (Here, $\tilde{S}$ is an
  $n$-fold covering of $S$.) It is straightforward to check that the maps $\tilde{u}$
  and $\tilde{\phi}$ again tautologically correspond.

  Now, it follows from the Atiyah-Singer index theorem that the index of the
  $\overline{\partial}$-operator at $\tilde{u}$ is $n$ times the index of the
  $\overline{\partial}$-operator at $u$. Alternately, this is immediate from the
  index formula on page 1018 of~\cite{L}, $\ind(u)=(g+\ell-1)/2-\chi(S)+2e(D(u)).$
\end{proof}

\emph{Remark.} Since the first version of this paper, S. Sarkar has proved a
combinatorial formula for the index of triangles (\cite{sarkarTriangles}). It should
be possible to use his formula to prove part~(\ref{item:ind}) of
Lemma~\ref{Lemma:CovTrip}.  We leave this to the interested reader.

\subsection{Bordisms and Heegaard triples} 
Recall that any bordism $W^4$ from $Y_1^3$ to $Y_2^3$ can be decomposed as a collection of \mbox{$1$-handle} attachments, followed by $2$-handle attachments, followed by $3$-handle attachments.  (See~\cite{GompfStipsicz} for an efficient exposition of handle decompositions and Kirby calculus.)  Attaching a \mbox{$1$-handle} to $W$ has the effect of either changing $\bdy W$ to $\bdy W\#(S^1\times S^2)$ or connect summing two connected components of $\bdy W$.  Attaching $2$-handles to $W$ has the effect on $\partial W$ of doing framed surgery along the attaching circles of the $2$-handles. Attaching a $3$-handle to $W$ has the opposite effect of attaching a $1$-handle: it either removes an $S^1\times S^2$ summand or it disconnects a connected sum.

In~\cite{holotri}, Ozsv\'ath and Szab\'o associate maps to $1$-handle, $2$-handle, and $3$-handle attachings, and show that the composition of these maps depends only on the bordism. The absolute $\qq$-grading, however, is defined using bordisms composed entirely of $2$-handles.  We will restrict our attention here to such bordisms, which we call for convenience \emph{link surgery bordisms}.

Fix a framed $n$-component link $\mathbb{L}=\{L_i\}$ in $Y$. By a \emph{bouquet for $\mathbb{L}$} we mean the union of $\mathbb{L}$ and a path from each component $L_i$ of $\mathbb{L}$ to a fixed reference point. Fix a bouquet $B(\mathbb{L})$ for $\mathbb{L}$.  A neighborhood $V$ of $B(\mathbb{L})$ is a handlebody of genus $n$. In general, $Y\smallsetminus V$ will not be a handlebody; however, it is possible to annex tubes to $V$ to produce a new handlebody $V'$ so that $Y\smallsetminus V'$ is a handlebody. Suppose $\bdy V'$ has genus $g$. We can choose $V'$ so that a small neighborhood of $L_i$ intersects $\bdy V'$ in a punctured torus $F_i$, with $F_i\cap F_j=\emptyset$ for $i\neq j$.

A Heegaard triple $(\Sigma,\alphas,\betas,\gammas,z)$ is \emph{subordinate to $\mathbb{L}$} if, for some choice of bouquet and $V'$ as above, there is an identification of $\Sigma$ with $\bdy V'$ such that
\begin{itemize}
\item Each $\alpha_i$ bounds a disk in $Y\smallsetminus V'$.
\item Each $\beta_i$ bounds a disk in $V'$.
\item For $1\leq i\leq n$, $\beta_i$ lies in $F_i$ and is a meridian of $L_i$.
\item For $n<i\leq g$, $\beta_i$ is disjoint from each $F_j$.
\item For $1\leq i\leq n$, $\gamma_i$ lies in $F_i$ and the homology class of $\gamma_i$ corresponds to the framing of $L_i$.
\item For $n<i\leq g$, $\gamma_i$ is a small perturbation of $\beta_i$.
\end{itemize}
This definition comes from~\cite[Section 4.1]{holotri}.

Observe that for such a Heegaard triple, $(\Sigma,\alphas,\beta_{n+1},\ldots,\beta_g)$ and $(\Sigma,\alphas,\gamma_{n+1},\ldots,\gamma_g)$ both specify $Y\smallsetminus \mathbb{L}=Y(\mathbb{L})\smallsetminus \mathbb{L}$, where $Y(\mathbb{L})$ denotes the $3$-manifold obtained by surgery along $\mathbb{L}$. Filling in the boundary of $Y\smallsetminus\mathbb{L}$ according to the $\beta_i$, $i\leq n$, gives back $Y$; filling in the boundary according to the $\gamma_i$, $i\leq n$, gives $Y(\mathbb{L})$.  That is, $Y_{\alpha,\beta}=Y$ and $Y_{\alpha,\gamma}=Y(\mathbb{L})$.

For any $B(\mathbb{L})$, there always exists a pointed Heegaard triple subordinate to it.  Note that if we start with a link surgery bordism $W$ and find a Heegaard triple $\TT_{\alpha,\beta,\gamma}$ subordinate to the corresponding link, then $W$ is obtained from $W_{\alpha,\beta,\gamma}$ by filling in $Y_{\beta,\gamma}=\#^{g-n}(S^1\times S^2)$ with $\natural^{g-n}(S^1\times\mathbb{D}^3)$.

\subsection{The absolute $\qq$-grading}\label{sec:AbsQGrading}
We will now define Ozsv\'{a}th and Szab\'{o}'s absolute $\qq$-grading, $\widetilde{\gr}$, on $\widehat{HF}(Y,\text{torsion})$, following the treatment in \cite[Section 7]{holotri}.  
First, for $\tt_i\in\spinc(Y_i)$, we define a \emph{$\spinc$ bordism from $(Y_1,\tt_1)$ to $(Y_2,\tt_2)$} to be a pair $(W,\ss)$ such that $W$ is a bordism from $Y_1$ to $Y_2$, and $\ss\in\spinc(W)$ satisfies $\ss|_{Y_i}=\tt_i$.  (We distinguish this concept from the concept of \emph{stable} $\spinc$ bordism, which will be relevant in Section 4.)
For any $\tt\in\spinc_\mathrm{tor}(Y)$, there exists a link surgery $\spinc$ bordism $(W,\ss)$ from $(S^3,\ss_0)$ to  $(Y,\tt)$, where $\ss_0$ is the unique $\spinc$ structure on $S^3$.  Let $\TT_{\alpha,\beta,\gamma}$ be a pointed Heegaard triple realizing $\ss$ and subordinate to the link.  Then $Y_{\alpha,\beta}=S^3$, $Y_{\beta,\gamma}=\#^k(S^1\times S^2)$ for some $k$, and $Y_{\alpha,\gamma}=Y$.  Choose $\xx_0\in\torus_{\alpha}\cap\torus_{\beta}$ and $\theta\in \torus_{\beta}\cap\torus_{\gamma}$ so that they have the same grading as the highest graded elements of $\widehat{HF}(S^3,\ss_0)$ and $\widehat{HF}(\#^k(S^1\times S^2)  ,\ss_0)$, respectively.  We say that $\xx_0$ and $\theta$ lie \emph{in the canonical degree}.
\begin{defin}
With notation as above, define $\tgr$ on $\widehat{CF}(\SS_{\alpha,\gamma},\tt)$ by
$$\widetilde{\gr}(\yy)=-\gr^\ss(\xx_0,\theta,\yy) + {c_{1}(\ss)^2-2\chi(W)-3\sigma(W)\over 4}$$
for $\yy\in\torus_{\alpha}\cap\torus_{\gamma}$ with $s(\yy)=\tt$.

\end{defin}
In \cite{holotri}, Ozsv\'{a}th and Szab\'{o} proved that this definition gives a well-defined absolute \mbox{$\qq$-grading} on $\widehat{HF}(Y,\text{torsion})$.  Observe that $\tgr$ is not obviously defined for a general pointed Heegaard diagram, but only those $\SS_{\alpha,\gamma}$ that arise from the above construction.  Call such a diagram a \emph{$\tgr$-admissible pointed Heegaard diagram for $Y$.}

For our purposes, we need to be able to work with $\tgr$ at the chain level.  The following cumbersome lemma may be thought of as a generalization of the definition above.  This argument is essentially the same as in the proofs of~\cite[Proposition 4.9, Lemma 7.5, and Theorem 7.1]{holotri}.  Since we want a chain-level statement, rather than a homology-level statement, we need to be slightly more precise.
\begin{lem}\label{tgrsurgery}
Let $(W,\ss)$ be a link surgery $\spinc$ bordism from $(Y_1,\tt_1)$ to $(Y_2,\tt_2)$.  Then there is a pointed Heegaard triple $\TT_{\alpha,\beta,\gamma}$ subordinate to a link inducing $W$, realizing $\ss$, such that $\SS_{\alpha,\beta}$ and  $\SS_{\alpha,\gamma}$ are $\tgr$-admissible for $Y_1=Y_{\alpha,\beta}$ and $Y_2=Y_{\alpha,\gamma}$, respectively. 

In this case, for all $\xx\in\torus_{\alpha}\cap\torus_{\beta}$ with $s(\xx)=\tt_1$ and $\yy\in\torus_{\alpha}\cap\torus_{\gamma}$ with $s(\yy)=\tt_2$,
\begin{equation}\label{tgrformula}
\widetilde{\gr}(\yy)=\widetilde{\gr}(\xx)-\gr^\ss(\xx,\theta,\yy)  + {c_{1}(\ss)^2-2\chi(W)-3\sigma(W)\over 4}.
\end{equation}
As before, $\theta\in \torus_{\beta}\cap\torus_{\gamma}$ lies in the canonical degree of $\widehat{HF}(Y_{\beta,\gamma}=\#^k(S^1\times S^2),\ss_0)$ for some~$k$.
\end{lem}
(We abuse notation here and elsewhere by identifying $\ss$ with $\ss|_{W_{\alpha,\beta,\gamma}}$.) 
\begin{proof}
Suppose that $\mathbb{L}'$ is a framed link in $S^3$ so that surgery on $\mathbb{L}'$ produces $Y_1$ and $\tt_1$ extends over the induced bordism from $S^3$ to $Y_1$.  Let $\mathbb{L}$ be any framed link in $Y_1$ inducing the link surgery bordism $W$.  By a small perturbation we can choose $\mathbb{L}$ to be disjoint from image of $\mathbb{L}'$ in $Y_1$, so that $\mathbb{L}$ is the image of some framed link in $S^3$.  Let $\mathbb{L}\mathbb{L}'$ be the union of $\mathbb{L}'$ and the preimage, in $S^3$, of $\mathbb{L}$.

Let $(\Sigma,\alphas,\deltas,\gammas)$ be a genus $g$ Heegaard triple subordinate to $\mathbb{L}\mathbb{L}'$, with the $\delta$ circles ordered so that $\delta_1,\ldots,\delta_m$ are meridians for components of $\mathbb{L}'$ and $\delta_{m+1},\ldots,\delta_n$ are meridians for components of the preimage of $\mathbb{L}$.  Let $\beta_i$ be a small isotopic translate of $\gamma_i$ for $1\leq i\leq m$;  let $\beta_i$ be a small isotopic translate of $\delta_i$ for $m< i\leq n$; and let $\beta_i$ be a small isotopic translate of both $\delta_i$ and $\gamma_i$ for $n< i\leq g$, all chosen so that each $\beta_i$ is transverse to $\alphas\cup\deltas\cup\gammas$.  Set $\betas=\beta_1\cup\cdots\cup\beta_g$. Then $(\Sigma,\alphas,\deltas,\betas)$ is a Heegaard triple subordinate to $\mathbb{L}'$, and $\TT_{\alpha,\beta,\gamma}=(\Sigma,\alphas,\betas,\gammas)$ is a Heegaard triple subordinate to $\mathbb{L}$.  We will now verify that $\TT_{\alpha,\beta,\gamma}$ satisfies the requirements in the statement of the lemma.

Let $W'$ denote the link surgery bordism induced by $\mathbb{L}'$, and $WW'$ the link surgery bordism induced by $\mathbb{L}\mathbb{L}'$, so that $WW'$ is the result of gluing $W$ to $W'$ along $Y_1$. Let $\ss'$ be a $\spinc$ structure on $W'$ extending $\tt_1$, and $\ss\ss'$ the $\spinc$ structure on $WW'$ induced by $\ss$ and $\ss'$.  (Note that it easy to make choices above so that $\ss$, $\ss'$, and $\ss\ss'$ are realized by the corresponding Heegaard triples.)

It is classical that 
\[
c_{1}(\ss\ss')^2-2\chi(WW')-3\sigma(WW')=c_1(\ss)^2+c_1(\ss')^2-2\chi(W)-2\chi(W')-3\sigma(W)-3\sigma(W').
\]
(In particular, for additivity of the signature, see~\cite[Section 7.1]{AtiyahSinger3}.)
To prove the result, then, it remains to check that the $\gr$-terms add.

Note that $Y_{\beta,\gamma}=\#^{g-n+m}(S^1\times S^2)$, $Y_{\delta,\beta}=\#^{g-m}(S^1\times S^2)$, and $Y_{\delta,\gamma}=\#^{g-n}(S^1\times S^2)$.  Choose $\theta\in\torus_{\beta}\cap\torus_{\gamma}$, $\theta'\in\torus_{\delta}\cap\torus_{\beta}$, and $\Theta\in\torus_{\delta}\cap\torus_{\gamma}$ lying in the canonical degree.

It is easy to check directly that in $(\Sigma,\deltas,\betas,\gammas)$, $\gr(\theta',\theta,\Theta)=0$. Finally, additivity properties of the index imply that for $\xx\in\torus_{\alpha}\cap\torus_{\beta}$ and $\yy\in\torus_{\alpha}\cap\torus_{\gamma}$ with $s(\xx)=\tt_1$ and $s(\yy)=\tt_2$, and $\xx_0\in \torus_{\alpha}\cap\torus_{\delta} $ in the canonical degree, 
\[
\gr(\xx_0,\theta\theta',\yy)+\gr(\theta',\theta,\Theta)=\gr(\xx_0,\theta',\xx)+\gr(\xx,\theta,\yy).
\]
where the suppressed superscripts are understood.  The result follows. 
\end{proof}

It is convenient to also understand how $\tgr$ behaves under connected 
sums. \begin{lem}\label{tgrsum} Let $\SS_1$ and $\SS_2$ be 
$\tgr$-admissible pointed Heegaard diagrams for $Y_1$ and $Y_2$.  Then 
$\SS_1\#\SS_2$ is a $\tgr$-admissible pointed Heegaard diagram for 
$Y_1\#Y_2$, and $$\widehat{CF}(\SS_1\# 
\SS_2)=\widehat{CF}(\SS_1)\otimes\widehat{CF}(\SS_2).$$ Furthermore, for 
any homogeneous elements $\xi_1\in\widehat{CF}(\SS_1)$ and 
$\xi_2\in\widehat{CF}(\SS_2)$ with $s(\xi_1)$ and $s(\xi_2)$ torsion, 
$$\tgr(\xi_1\otimes\xi_2)=\tgr(\xi_1)+\tgr(\xi_2).$$ \end{lem}

This lemma can be proved by considering a link surgery bordism from $Y_1$ to
$Y_1\#Y_2$ and applying the same reasoning used in the proof of the previous
lemma. Here, it is convenient to choose the $\mathbb{L}$ and $\mathbb{L}'$ so that
$\mathbb{LL'}\subset S^3$ is a split link. Then the corresponding Heegaard triple
subordinate to $\mathbb{LL'}$ can be chosen to be a connected sum of a triple for
$Y_1$ and a triple for $Y_2$.

Finally, given a pointed Heegaard diagram $\SS$ for $Y$, observe that if 
$\tt\in\spinc_{\mathrm{tor}}(Y)$ and $\widehat{HF}(Y,\tt)\neq0$, then 
$\tgr$ on $\widehat{HF}(Y,\tt)$ determines $\tgr$ on 
$\widehat{CF}(\SS,\tt)$.

\section{Proof that the two constructions agree}

The purpose of this section is to prove the following theorem.
\begin{thm}\label{mainagree}
$\Gr=\widetilde{\gr}$ as relative $\qq$-gradings on $\widehat{HF}(Y,\mathrm{torsion})$.
\end{thm}

The idea of the proof is the following.  The grading $\tgr$ is well-behaved under $\spinc$ bordisms.  In particular, if we know $\tgr$ on $\widehat{CF}(Y,\tt)$, then we can determine $\tgr$ on any pair that is $\spinc$ bordant to $(Y,\tt)$.  The definition of $\tgr$ in Section 3.3 works because there is only one $\spinc$ bordism class in dimension three.  The relative grading $\Gr$ has similar good behavior under $\spinc$ bordism, but only those bordisms admitting $\zz/n$-covers.  In particular if we know $\Gr$ on $\widehat{CF}(Y,\tt)\oplus\widehat{CF}(Y,\tt')$, then we can determine $\Gr$ on something that is ``suitably'' $\spinc$ bordant to $(Y,\tt,\tt')$.  Therefore, if we can show that $\Gr$ is consistent with $\tgr$ on a representative of each ``suitable'' bordism class (it turns out that lens spaces will suffice), then the result will follow.

We begin by showing that $\Gr$ has the desired behavior under connected sums.
\begin{lem}\label{Grsum}
Let $\SS_1$ and $\SS_2$ be pointed Heegaard diagrams for $Y_1$ and $Y_2$.  Then $\SS_1\#\SS_2$ is a pointed Heegaard diagram for $Y_1\#Y_2$, and
$$\widehat{CF}(\SS_1\# \SS_2)=\widehat{CF}(\SS_1)\otimes\widehat{CF}(\SS_2).$$
Furthermore, for any homogeneous elements $\xi_1,\xi_1'\in\widehat{CF}(\SS_1)$ and $\xi_2,\xi_2'\in\widehat{CF}(\SS_2)$ with $s(\xi_1)$, $s(\xi_1')$, $s(\xi_2)$, and $s(\xi_2')$ torsion,
$$\Gr(\xi_1\otimes\xi_2,\xi_1'\otimes\xi_2')=\Gr(\xi_1,\xi_1')+\Gr(\xi_2,\xi_2').$$ 
\end{lem}
\begin{proof}
Consider a cover of $Y_1$ for which $\ss(\tilde{\xi_1})=\ss(\tilde{\xi_1'})$ and a cover of $Y_2$ for which $\ss(\tilde{\xi_2})=\ss(\tilde{\xi_2'})$.  From these two covers, we can construct a cover of $Y_1\#Y_2$ for which the computation becomes obvious.
\end{proof}
Combining this lemma with Lemma \ref{tgrsum}, we have the following.
\begin{cor}\label{sumagree}
If $\Gr=\widetilde{\gr}$ as relative $\qq$-gradings on $\widehat{HF}(Y_1,\mathrm{torsion})$ and $\widehat{HF}(Y_2,\mathrm{torsion})$, then they are equal on $\widehat{HF}(Y_1\# Y_2,\mathrm{torsion})$.
\end{cor}

Now observe that trivially, $\Gr=\tgr$ as relative $\qq$-gradings on $\widehat{HF}(S^1\times S^2,\ss_0)\isom \zz\oplus\zz$.  Therefore Lemmas \ref{tgrsum} and \ref{Grsum} give us the following. 
\begin{cor}\label{handleagree}
Let $\tt,\tt'\in\spinc_\mathrm{tor}(Y)$.  Then $\Gr=\widetilde{\gr}$ as relative $\qq$-gradings on $\widehat{HF}(Y,\tt)\oplus\widehat{HF}(Y,\tt')$ if and only if they are equal on $\widehat{HF}(Y\#(S^{1}\times S^{2}),\tt\#\ss_{0})\oplus\widehat{HF}(Y\#(S^{1}\times S^{2}),\tt'\#\ss_{0})$. 
\end{cor}

The following lemma will be useful for producing equivariant $\spinc$ bordisms. When $\tt$ is a $\spin$ structure we will denote the induced $\spinc$ structure by $\tt$ as well. 

\begin{lem}\label{Lem:Bordism}
Given $f:Y\to B\zz/n$ and a $\spin$ structure $\tt$ on $Y$, there exists a bordism $F:W\to B\zz/n$ from $f$ to a disjoint union of maps $f_1:L(n,1)\to B\zz/n$, and there exists $\ss\in\spinc(W)$ such that $\ss|_Y=\tt$.
\end{lem}
\begin{proof}
First, recall the following (stable) $\spin$ bordism groups of a point.
\begin{eqnarray*}
\Omega_0^{\spin}&=&\zz\\ 
\Omega_1^{\spin}&=&\zz/2\\
\Omega_2^{\spin}&=&\zz/2\\
\Omega_3^{\spin}&=&0.
\end{eqnarray*}
We use the Atiyah-Hirzebruch spectral sequence in order to understand the third (stable) $\spin$ bordism group of $B\zz/n$, $\Omega_3^\spin(B\zz/n)$.  We know that each class in $\Omega_3^\spin(B\zz/n)$ can be described by an element of

$$\bigoplus_{p+q=3} E^2_{p,q}\\
\isom H_0(B\zz/n; \Omega_3^{\spin})\oplus H_1(B\zz/n; \Omega_2^{\spin})\oplus
H_2(B\zz/n; \Omega_1^{\spin})\oplus H_3(B\zz/n;\Omega_0^{\spin}).
$$

Consider the last summand  $H_3(B\zz/n; \Omega_0^{\spin})\isom \zz/n$. 
Choose a model of $B\zz/n$ with $L(n,1)$ as its 3-skeleton; this map $f_1:L(n,1)\to B\zz/n$, together with a $\spin$ structure on $L(n,1)$ gives us a generator of $H_3(B\zz/n; \Omega_0^{\spin})$.

We now consider the other three summands. The first is obviously zero, since $\Omega_3^{\spin}=0$. Since $\Omega_1^{\spinc}=0$ and $\Omega_2^{\spinc}=\zz$, the maps $\Omega_1^{\spin}\to \Omega_1^{\spinc}$ and $\Omega_2^{\spin}\to \Omega_2^{\spinc}$ are both zero. By  naturality of the Atiyah-Hirzebruch spectral sequence, it follows that any element of $H_1(B\zz/n; \Omega_2^{\spin})\oplus H_2(B\zz/n; \Omega_1^{\spin})$ is \emph{stably} $\spinc$ bordant to 0.

We would be done now, except that that stable $\spin$ bordisms and stable $\spinc$ bordisms are manifolds with $\spin$ or $\spinc$ structures on the stable normal bundle, not the tangent bundle. For $\spin$ bordism this is no additional complication: there is a one-to-one correspondence between $\spin$ structures on the tangent bundle and $\spin$ structures on the stable normal bundle (see \cite[Proposition 2.15]{LawsonMichelsohn}).  In particular, $(f,\tt)$ must be $\spin$ bordant to something in $\Omega_3^\spin(B\zz/n)$

For $\spinc$ structures on $3$-manifolds, there is no such correspondence. However, for manifolds of dimension less than three, there is such a correspondence: in these dimensions, $\spinc$ structures on the tangent bundle and $\spinc$ structures on the stable normal bundle are both classified by their first Chern classes.  Therefore any element of $H_1(B\zz/n; \Omega_2^{\spin})\oplus H_2(B\zz/n; \Omega_1^{\spin})$ is honestly $\spinc$ bordant to 0.  The result follows.

\end{proof}

\begin{cor}\label{Cor:EquivariantBordism}Fix a $3$-manifold $Y$, a $\spin$ structure $\tt$ on $Y$, and a torsion $\spinc$ structure $\tt'$ on $Y$. Suppose that $\langle \tt-\tt'\rangle=\zz/n\subset H^2(Y;\zz)$. Let $\tilde{Y}\to Y$ be a $\zz/n$-covering space as in the statement of Corollary \ref{lift}.  Then for some $\uu,\uu'\in\spinc(\coprod^m(L(n,1)))$, there are $\spinc$ bordisms $(W,\ss)$ from $(\coprod^m(L(n,1)),\uu)$ to $(Y,\tt)$ and $(W,\ss')$ from $(\coprod^m(L(n,1)),\uu')$ to $(Y,\tt')$, and another $\spinc$ bordism $(\tilde{W},\tilde{\ss})$ from $S^3$ to $\tilde{Y}$ covering $(W,\ss)$ and $(W,\ss')$.  That is,
\protect{
\[
\xymatrix{
(\coprod^m(L(n,1)),\uu) \ar@{-}[rr]|(.57){(W,\ss)} && (Y,\tt)\\
(S^3,{\ss_0}) \ar@{-}[rr]|(.57){(\tilde{W},\tilde{\ss})} \ar[d]\ar[u] &\ar[u]\ar[d]& (\tilde{Y},\tilde{\tt})\ar[d]\ar[u]\\
(\coprod^m(L(n,1)),\uu') \ar@{-}[rr]|(.57){(W,\ss')} && (Y,\tt')}
\]}
Further, $\ss-\ss'$ is a torsion element of $H^2(W;\zz)$.
\end{cor}

\begin{proof}
  Let $\tilde{Y}\to Y$ be the covering map given by Corollary \ref{lift}, and let
  $f:Y\to B\zz/n$ be its classifying map.  Then by Lemma~\ref{Lem:Bordism} there is a
  $\spinc$-bordism $F: (W,\ss)\to B\zz/n$ from $f: (Y,\tt)\to B\zz/n$ to a
  disjoint union of $m$ maps of the form $f_1: (L(n,1),\uu)\to B\zz/n$. (Here,
  $\uu$ is some fixed $\spin$-structure on $L(n,1)$.) Recall from the proof of Lemma
  \ref{topology} that $\tt-\tt'=f^*b$ for some $b\in H^2(B\zz/n;\zz)$.  Setting
  $\ss'=\ss-F^*b$, we see that $\ss'|_Y=\tt'$.  Setting $\tilde{W}=F^*(E\zz/n)$, we
  see that $\ss$ and $\ss'$ pullback to the same $\tilde{\ss}\in\spinc(\tilde{W})$.
\end{proof}

As mentioned earlier, we will need to know that the two gradings agree for lens spaces.
\begin{lem}\label{LensAgree}
$\Gr=\widetilde{\gr}$ as relative $\qq$-gradings on $\widehat{HF}(L(n,1))$.
\end{lem}
We defer the proof of this lemma until Section 5, where it follows from Proposition~\ref{LensOurs} and Corollary~\ref{LensOS}.

\begin{proof}[Proof of Theorem \ref{mainagree}]
For the sake of clarity, assume for now that there exists a $\spin$ structure $\tt$ on $Y$ such that
$\widehat{HF}(Y,\tt)\neq 0$; this is necessarily the case, for instance, if $Y$ is a rational homology sphere.  At the end of the proof we will explain how to remove this assumption.  With this assumption, it is sufficient to show that for all $\tt'\in\spinc_\mathrm{tor}(Y)$, $\Gr=\widetilde{\gr}$ as relative $\qq$-gradings on $\widehat{HF}(Y,\tt)\oplus\widehat{HF}(Y,\tt')$.

Let $W$, $\tilde{W}$, $\ss$, $\ss'$, $\uu$, $\uu'$, and $m$ be as in the statement of Corollary \ref{Cor:EquivariantBordism}.  As described in Section~3.2, we can decompose $W$ into a bordism from $\coprod^m(L(n,1))$ to $Y_{1}=(\#^m(L(n,1)))\#(\#^{k_1}(S^1\times S^2))$,  a link surgery bordism $W_1$ from $Y_{1}$ to $Y_{2}=Y\#(\#^{k_2}(S^1\times S^2))$, and a bordism from $Y_2$ to $Y$.  By the previous lemma, we know that $\Gr=\widetilde{\gr}$ as relative $\qq$-gradings on $\widehat{HF}(L(n,1))$.  Then by Corollaries \ref{sumagree} and \ref{handleagree}, we know that $\Gr=\widetilde{\gr}$ as relative $\qq$-gradings on $\widehat{HF}(Y_{1},\text{torsion})$.

Now choose a pointed Heegaard triple $\TT_{\alpha,\beta,\gamma}$ corresponding to $W_1$ as in the statement of Lemma \ref{tgrsurgery} such that $\TT_{\alpha,\beta,\gamma}$ realizes both $\ss|_{W_1}$ and $\ss'|_{W_1}$.  Let $\xx,\xx'\in\torus_{\alpha}\cap\torus_{\beta}$ with $s(\xx)=(\#\uu)\#\ss_0$ and  $s(\xx')=(\#\uu')\#\ss_0$, and let $\yy,\yy'\in\torus_{\alpha}\cap\torus_{\gamma}$ with $s(\yy)=\tt\#\ss_0$ and $s(\yy')=\tt'\#\ss_0$.  Then
\begin{eqnarray}
\widetilde{\gr}(\yy)-\widetilde{\gr}(\yy')&=&\widetilde{\gr}(\xx)-\widetilde{\gr}(\xx')-(\gr^\ss(\xx,\theta,\yy)-\gr^{\ss'}(\xx',\theta,\yy'))+ {c_{1}(\ss|_{W_1})^2-c_{1}(\ss'|_{W_1})^2\over 4}\nonumber\\
&=&\widetilde{\gr}(\xx)-\widetilde{\gr}(\xx')-\gr^\ss(\xx,\theta,\yy)+\gr^{\ss'}(\xx',\theta,\yy')
\label{eqn:2}
\end{eqnarray}
where the last line follows because $\ss-\ss'$ is torsion.
\begin{claim}
$\Gr(\yy,\yy')=\Gr(\xx,\xx')-\gr^\ss(\xx,\theta,\yy)+\gr^{\ss'}(\xx',\theta,\yy')$.
\end{claim}
Using the cover $\tilde{W}\to W$, construct a covering Heegaard triple $\tilde{\TT}_{\tilde{\alpha},\tilde{\beta},\tilde{\gamma}}$ and compute:
\begin{eqnarray*}
\Gr(\yy,\yy')&=&{1\over n}\gr(\tilde{\yy},\tilde{\yy}')\\
&=&{1\over n}(\gr(\tilde{\xx},\tilde{\xx}')-\gr^{\pi^*\ss}(\tilde{\xx},\tilde{\theta},\tilde{\yy})+\gr^{\pi^*\ss'}(\tilde{\xx}',\tilde{\theta},\tilde{\yy}'))       \text{ by additivity of $\gr$}\\
&=&\Gr(\xx,\xx')-\gr^\ss(\xx,\theta,\yy)+\gr^{\ss'}(\xx',\theta,\yy')\text{ by Lemma~\ref{Lemma:CovTrip}.}
\end{eqnarray*}
Since $\widehat{HF}(Y_1,(\# \uu)\#\ss_0)$ and $\widehat{HF}(Y_1,(\# \uu')\#\ss_0)$ are nonzero and $\Gr=\widetilde{\gr}$ as relative \mbox{$\qq$-gradings} on $\widehat{HF}(Y_{1},\text{torsion})$, it follows that $\Gr(\xx,\xx')=\widetilde{\gr}(\xx')-\widetilde{\gr}(\xx)$.  Equation~(\ref{eqn:2}) and the claim now show that $\Gr=\widetilde{\gr}$  as relative $\qq$-gradings on $\widehat{CF}(\SS_{\alpha,\gamma},\tt\#\ss_0)\oplus\widehat{CF}(\SS_{\alpha,\gamma},\tt'\#\ss_0)$.  Consequently they are equal on $\widehat{HF}(Y_{2},\tt\#\ss_0)\oplus\widehat{HF}(Y_{2},\tt'\#\ss_0)$, and
Corollary \ref{handleagree} then shows that they are equal on  $\widehat{HF}(Y,\tt)\oplus 
\widehat{HF}(Y,\tt')$.  

We now deal with the general case when the simplifying assumption fails.  Choose a
$\spin$ structure $\tt$ on $Y$.  To prove the theorem, it is sufficient to show that
for all $\tt',\tt^\circ\in\spinc_\mathrm{tor}(Y)$, $\Gr=\widetilde{\gr}$ as relative
$\qq$-gradings on $\widehat{HF}(Y,\tt')\oplus\widehat{HF}(Y,\tt^\circ)$.  We will do
this by constructing a Heegaard diagram $\SS_{\alpha,\gamma}$ for
$Y\#(\#^{k_2}S^1\times S^2)$ for which we know the relative gradings $\Gr$ and $\tgr$
agree on
$\widehat{CF}(\SS_{\alpha,\gamma},\tt\#\ss_0)\oplus\widehat{CF}(\SS_{\alpha,\gamma},\tt'\#\ss_0)$
and also agree on
$\widehat{CF}(\SS_{\alpha,\gamma},\tt\#\ss_0)\oplus\widehat{CF}(\SS_{\alpha,\gamma},\tt^\circ\#\ss_0)$. By
Corollary~\ref{handleagree}, this is sufficient to prove that $\Gr$ and $\tgr$ agree
as relative gradings on $\widehat{HF}(Y,\tt')\oplus\widehat{HF}(Y,\tt^\circ)$.

Let $W$, $\tilde{W}$, $\ss$, $\ss'$, $\uu$, $\uu'$, $m$, $Y_1$ and $Y_2$ be the
objects constructed before using $\tt$ and $\tt'$. Let $W^\circ$, $\tilde{W}^\circ$,
$\ss^\circ$, $\ss'^\circ$, $\uu^\circ$, $(\uu^\circ)'$, $m^\circ$, $Y_1^\circ$ and
$Y_2^\circ$ be the corresponding objects constructed using $\tt$ and $\tt^\circ$. By
adding some canceling $1$- and $3$-handles to the decompositions of $W$ or $W^\circ$
if necessary, we may assume that $Y_2=Y_2^\circ$.

Now, choose Heegaard triples $\TT_{\alpha,\beta,\gamma}$ and
$\TT^\circ_{\alpha^\circ,\beta^\circ,\gamma^\circ}$ for $W_2$ and $W_2^\circ$
respectively, with the property that
$\SS_{\alpha,\gamma}=\SS^\circ_{\alpha^\circ,\gamma^\circ}$. This is possible by the argument
used in Lemma~\ref{tgrsurgery}. That is, we view $W_2$ and $W_2^\circ$ as given by
surgery on disjoint links $\mathbb{L}$ and $\mathbb{L'}$ in $Y_2$. Then it is easy to construct the desired
triple diagrams from a diagram subordinate to $\mathbb{L}\cup\mathbb{L'}$. Further,
arrange that $\widehat{CF}(\SS_{\alpha,\gamma},\tt\#\ss_0)$ is nonzero, \textit{i.e.}, that
there is some generator in the $\spinc$-structure $\tt\#\ss_0.$

As before, it follows from Corollaries~\ref{sumagree} and~\ref{handleagree} that
$\Gr=\tgr$ as relative $\qq$-gradings on both
$\widehat{CF}(\SS_{\alpha,\beta},\tt\#\ss_0)\oplus\widehat{CF}(\SS_{\alpha,\beta},\tt'\#\ss_0)$
and
$\widehat{CF}(\SS^\circ_{\alpha,\beta^\circ},\tt\#\ss_0)\oplus\widehat{CF}(\SS^\circ_{\alpha,\beta^\circ},\tt^\circ\#\ss_0)$.
Then, by Formula~(\ref{eqn:2}) and the Claim above, it follows that $\Gr=\tgr$ as relative $\qq$-gradings $\widehat{CF}(\SS_{\alpha,\gamma},\tt\#\ss_0)\oplus\widehat{CF}(\SS_{\alpha,\gamma},\tt'\#\ss_0)$
and
$\widehat{CF}(\SS_{\alpha,\gamma},\tt\#\ss_0)\oplus\widehat{CF}(\SS_{\alpha,\gamma},\tt^\circ\#\ss_0)$.
It then follows that $\Gr=\tgr$ as relative $\qq$-gradings on
$\widehat{HF}(\SS_{\alpha,\gamma},\tt'\#\ss_0)\oplus\widehat{HF}(\SS_{\alpha,\gamma},\tt^\circ\#\ss_0)$.
The result now follows from Corollary~\ref{handleagree}.
\end{proof}

As discussed in Section~\ref{sec:GradingReview}, there is a group $\widehat{HF}(Y,\ell)$ that is independent of choice of $\ell$-pointed Heegaard diagram.  The invariant absolute grading $\tgr$ was originally only defined on $\widehat{HF}(Y,1)$ \cite{holotri}.   In light of the two main theorems of this paper, Theorems~\ref{main} and~\ref{mainagree}, one might hope that some generalization of $\tgr$ to an invariant absolute grading of $\widehat{HF}(Y,\ell)$ would behave nicely with respect to coverings.  Specifically, given a one-pointed Heegaard diagram $\SS$ for $Y$ and an $n$-fold cover of $Y$, one might hope that $\tgr(\xx)=\frac{1}{n}\tgr(\tilde{\xx})$ (whenever $\xx$ and $\tilde{\xx}$ are both nontrivial in homology).  Unfortunately, the example of lens spaces shows that this is impossible.

Explicitly, consider the standard one-pointed Heegaard diagrams $\SS_1$ and $\SS_2$ for $L(5,1)$ and $L(5,2)$, respectively, together with the covering maps from $S^3$ to $L(5,1)$ and $L(5,2)$.  One can find $\xx_1\in\widehat{CF}(\SS_1)$ and $\xx_2\in\widehat{CF}(\SS_2)$ with the properties that:
\begin{itemize}
\item $\tilde{\xx}_1$ and $\tilde{\xx}_2$ define the same nontrivial element of
  $\widehat{HF}(S^3,5)$.
\item $\xx_1$ and $\xx_2$ are both nontrivial in homology, and $\tgr(\xx_1)\neq\tgr(\xx_2)$.
\end{itemize}

\section{Lens spaces and gradings}

\piccaption{The lens space $-L(5,1)$.\label{Figure:Lens}}
\parpic{
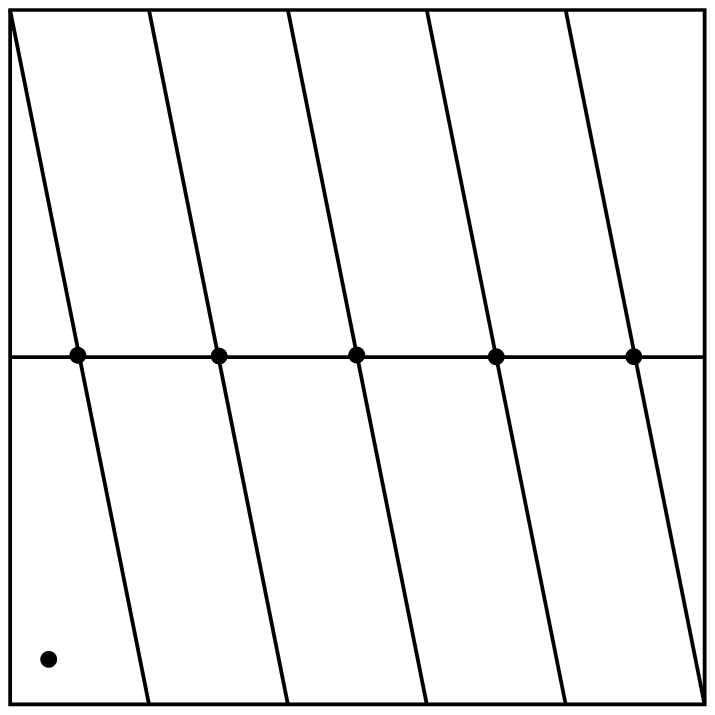}
Let $p$ and $q$ be relatively prime positive integers with $p>q$.  In this section we will compute $\Gr$ on $L(p,q)$ and show that the answer agrees with the relative \mbox{$\qq$-grading} induced by $\widetilde{\gr}$, as computed by Ozsv\'ath and Szab\'o in~\cite[Proposition 4.8]{absGrade}.  We first describe a pointed Heegaard diagram $\mathcal{S}=(\Sigma,\alpha,\beta,z)$ for $-L(p,q)$.  Let $\Sigma$ be a torus, realized as the square $[0,1]\times[0,1]$ with edges identified.  Let $\alpha$ be the horizontal circle $y=1/2$, and let $\beta$ be any circle with slope $-p/q$.  Choose $z$ to be any point in $\Sigma\smallsetminus\alpha\smallsetminus\beta$ lying below $y=1/2$. There are $p$ connected components of $\Sigma\smallsetminus\alpha\smallsetminus\beta$.  We label their closures $D_0,\ldots,D_{p-1}$, going from left to right, starting with $z\in D_0$.  (Note that this is not the same notational convention used in Section 2.)  There are $p$ elements of $\alpha\cap\beta$.  We label them $x_0,\ldots,x_{p-1}$, going from left to right along $\alpha$, such that $x_0$ is the upper right vertex of $D_0$.  See Figure 1.

It is easy to see that $\epsilon(\xx_i,\xx_j)\neq0$ if $i\neq j$, so each $\ss_i=s(x_i)$ is distinct, and these are all of the $\spinc$ structures on $-L(p,q)$.  Consequently, the differential on $\widehat{CF}(\mathcal{S})$ is trivial, and we can identify $\widehat{CF}(\mathcal{S})$ with $\widehat{HF}(-L(p,q))$.  Ozsv\'{a}th and Szab\'{o} proved the following inductive formula for their absolute \mbox{$\qq$-grading} on $\widehat{HF}(-L(p,q))$ in~\cite[Proposition 4.8]{absGrade}. Since their formula relates $\widetilde{\gr}$ on different lens spaces, we will use the notation $\widetilde{\gr}_{p,q}$ for the absolute \mbox{$\qq$-grading} on $\widehat{HF}(-L(p,q))$.
\begin{prop}\label{prop:theirformula}
In the setup described above, for $0\leq i< p+q$,
$$\widetilde{\gr}_{p,q}(x_{i\!\!\!\mod    {p}})=\left({ pq-(2i+1-p-q)^2\over 4pq}\right)-\widetilde{\gr}_{q,p\!\!\!\mod    {q}}(x_{i\!\!\!\mod    {q}}).$$
\end{prop}

We can use this formula to derive a non-inductive formula for the relative \mbox{$\qq$-grading} induced by $\tgr$.
\begin{cor}\label{LensOS}
Fix $p$ and $q$, and let $\widetilde{\gr}$ refer to the absolute \mbox{$\qq$-grading} on $\widehat{HF}(-L(p,q))$.  For 
$0\leq i<p$,
$$\widetilde{\gr}(x_{i+q\!\!\!\mod      {p}}) - \widetilde{\gr}(x_i)={1\over p}(p-1-2i).$$
\end{cor}
Note that this formula completely determines the relative \mbox{$\qq$-grading} on $\widehat{HF}(-L(p,q))$.
\begin{proof}
By the previous proposition,
\begin{eqnarray*}
\widetilde{\gr}_{p,q}(x_{i+q\!\!\!\mod      p})-\widetilde{\gr}_{p,q}(x_{i\!\!\!\mod      p})
&=&\left({ pq-(2(i+q)+1-p-q)^2\over 4pq}\right)-\widetilde{\gr}_{q,p\!\!\!\mod      q}(x_{i+q\!\!\!\mod      q})\\
& &-\left({ pq-(2i+1-p-q)^2\over 4pq}\right)+\widetilde{\gr}_{q,p\!\!\!\mod      q}(x_{i\!\!\!\mod      q})\\
&=&\left({ -(2i+1-p+q)^2\over 4pq}\right)+\left({(2i+1-p-q)^2\over 4pq}\right)\\
&=&{2(2i+1-p)(-2q)\over 4pq}\\
&=&{1\over p}(p-1-2i).
\end{eqnarray*}
\end{proof}
We now derive the same formula for $\Gr$. 
\begin{prop}\label{LensOurs}
Let $\Gr$ be the relative \mbox{$\qq$-grading} on $\widehat{HF}(-L(p,q))$ defined in Section 2.  For $0\leq i<p$,
$${\Gr}(x_{i+q\!\!\!\mod      p},x_i)={1\over p}(p-1-2i).$$
\end{prop}
\begin{proof}
Since $p$ and $q$ are now fixed, we can interpret the indices for $x$ and $D$ as being defined mod $p$.
In order to compute $\Gr(x_{i+q},x_i)$ using Proposition
 \ref{Ntheory}, we need an element of $\hat{\pi}_2(px_{i+q},px_{i})$. 
\begin{claim}
For $0\leq i<p$, let $A_i=\sum_{j=1}^{p} (p-(i+j))D_{i+j}$.  Then $A_i\in \hat{\pi}_2(px_{i+q},px_{i})$.
\end{claim}
To see that $A\in \hat{C}_2(\SS)$, note that the coefficient of $D_0=D_p$ is zero.  The proof that $A\in{\pi}_2(px_{i+q},px_{i})$ is a simple computation using the fact that 
$$\partial\partial_\alpha D_i=x_{i+q}-x_{i+q-1}+x_{i-1}-x_{i}.$$
We can now plug $A_i$ into the formula in Proposition \ref{Ntheory}.  Noting that $e(D_i)=0$, $x_i$ touches $D_{p+i}$, $D_{i+1}$, $D_{p+i-q}$, and $D_{p+i-q+1}$, and $x_{i+q}$ touches $D_{i+q}$, $D_{i+q+1}$, $D_{p+i}$, and $D_{i+1}$, we compute
\begin{eqnarray*}
\Gr(x_{i+q},x_i)&=&{1\over p}[e(A_i)+n_{x_{i+q}}(A_i)+n_{x_{i}}(A_i)]\\
&=&{1\over p}\left(0+{1\over4}[(p-i-q)+(p-i-q-1)+(-i)+(p-i-1)]\right.\\
& &+\left.{1\over 4}[(-i)+(p-i-1)+(-i+q)+(-i+q-1)]\right)\\
&=&{1\over p}(p-1-2i).
\end{eqnarray*}
\end{proof}

\section{Directions for future research}
This paper is an offshoot of an ongoing attempt by the authors to understand the relationship between the Floer homology of a space and the Heegaard Floer homology of its finite covering spaces, the main technical goal of which is a localization theorem. This is one direction for further study.

More directly related to gradings, it would be interesting to extend our definition to non-torsion $\spinc$ structures.  This has two meanings.  The simpler is that if $\tt_1-\tt_2\in H^2(Y)$ is torsion then one can find a covering space $p:\tilde{Y}\to Y$ such that $p^*\tt_1=p^*\tt_2$; one could try to use this covering space to define a relative grading between generators of $\widehat{CF}(Y,\tt_1)$ and $\widehat{CF}(Y,\tt_2)$. The difficulty arises from the fact that when $\tt$ is non-torsion, the relative grading on $\widehat{CF}(Y,\tt)$ is only defined modulo $\gcd_{A\in H_2(Y)}\{\langle c_1(\tt),A\rangle\}$.  Hopefully, this difficulty could either be overcome or exploited.

A larger generalization would be to obtain a relative grading between $\spinc$ structures with non-torsion difference by considering infinite covering spaces.  In light of~\cite{AtiyahPatodiSinger2}, it seems likely that here the real numbers, rather than the rationals, would come into play.  Conceivably, attempts to generalize the index could lead to some kind of ``$\ell^2$ Heegaard Floer homology'' for infinite covering spaces.

Finally, it would be interesting to explore how covering spaces might be used to construct the absolute \mbox{$\qq$-grading} $\tgr$, as well as the induced relative $\qq$-grading.  In view of the example of lens spaces, such a construction might require significant new ideas.

\bibliographystyle{alpha}
\bibliography{researchRL}

\end{document}